\documentclass{article}
\usepackage{latexsym,amssymb}
\usepackage[dvips]{graphicx}
\usepackage{subfigure}

\newcommand{\C}{\mathbb{C}}
\newcommand{\R}{\mathbb{R}}

\newcommand{\<}{\left<}
\newcommand{\m}{\right>}

\newcommand{\G}{\mathcal{N}}
\def\NN{\mathcal{N}}

\newcommand{\p}{\mathcal{P}}

\newcommand{\J}{J}
\newcommand{\Fix}{\mathop{\rm Fix}}
\newcommand{\Ker}{\mathop{\rm Ker}}

\newcommand{\Coad}{\mathop{\rm Coad}\nolimits}

\newcommand{\se}{\mathfrak{se}(2)}
\newcommand{\sodeux}{\mathfrak{so}(2)}
\newcommand{\g}{\mathfrak{g}}

\def\rum{relative equilibrium }
\def\rumb{relative equilibrium}
\def\ria{relative equilibria }
\def\riab{relative equilibria}
\def\rmk{\noindent\textbf{Remark. }}

\def\mom{momentum map }
\def\tet{\theta}
\def\al{\alpha}
\def\1N{1\dots N}
\def\la{\lambda}
\def\om{\omega}
\def\fer{\Omega}
\newtheorem{theo}{Theorem}[section]
\newtheorem{coro}[theo]{Corollary}
\newtheorem{prop}[theo]{Proposition}

\newenvironment{proof}{\addvspace\baselineskip\noindent{\textbf{Proof.}}\quad}{\hspace*{\fill} $\Box$\par\addvspace\baselineskip}

\newcommand{\bfig}{\begin{figure}[htbp]}
\newcommand{\efig}{\end{figure}}
\newcommand{\bc}{\begin{center}}
\newcommand{\ec}{\end{center}}
\newcommand{\beq}{\begin{equation}}
\newcommand{\eeq}{\end{equation}}
\newcommand{\pg}{\textquotedblleft}
\newcommand{\pd}{\textquotedblright}

\begin{document}

\bc
{\Large Point vortices on a rotating sphere}

\vspace{0.3cm}

{\large Fr\'{e}d\'{e}ric Laurent-Polz}\\
\textit{Institut Non Lin\'eaire de Nice, Universit\'e de Nice, France}\\
laurent@inln.cnrs.fr
\ec

\vspace{0.5cm}

\bc
\noindent\textbf{Abstract}\\
\ec
\noindent We study the dynamics of $N$ point vortices on a rotating sphere. The Hamiltonian system becomes
infinite dimensional due to the non-uniform background vorticity coming from the Coriolis force.
We prove that a \rum formed of latitudinal rings of identical vortices for the non-rotating
sphere persists to be a \rum when the sphere rotates. We study the nonlinear stability of a polygon
of planar point vortices on a rotating plane in order to approximate the corresponding \rum on the
rotating sphere when the ring is close to the pole.
We then perform the same study for geostrophic vortices.
To end, we compare our results to the observations on the southern hemisphere atmospheric
circulation.

\bigskip


\noindent AMS classification scheme number: 70E55, 70H14, 70H33\\

\noindent PACS classification scheme number: 45.20.Jj, 45.50.Jf, 47.20.Ky, 47.32.Cc


\section{Introduction}

The interest of studying point vortices on a rotating sphere is clearly geophysical.
This may permit also to understand the motion of concentrated regions of vorticity on the surface
of planets such as Jupiter \cite{DL93}.
The literature is now numerous on point vortices on a \emph{non}-rotating sphere
\cite{B77,KN98,KN99,KN00,N00,PM98,LMR00,BC01,LP02} but only few papers consider a rotating sphere
\cite{F75,B77,B85,KR89,DP98}.
In \cite{F75}, the interaction of three identical point vortices equally spaced on the same
latitude is investigated \emph{via} the $\beta$-plane approximation. It is shown --- under some
additional assumptions --- that this configuration is an equilibrium, and that its linear stability
depends on the strength of the vortices:
linearly stable for a negative or a strongly positive strength, linearly unstable otherwise.
In \cite{B77}, the equations of motion on a rotating sphere are given, while in \cite{B85} the
motion of a single point vortex is given.
It appears that a single vortex moves westward and northward as a hurricane does.
In \cite{KR89}, the approach is completely different from ours and \cite{B77};
they proved the existence of \ria formed of two or three vortices of possibly non-identical
vorticities.
In \cite{DP98}, they model the background vorticity coming from the rotation by
latitudinal strips of constant vorticities and they study the motion of a vortex pair
(point-vortex pair as well as patch-vortex pair). A vortex pair is a solution formed of two vortices
with opposite sign vorticities rotating around the North pole.
A vortex pair moving eastward or strongly westward is stable, unless the vortex sizes are too
large. Different types of instabilities are described for weak westward pairs.

In this paper, we first recall some basics notions of geometric mechanics in Section \ref{geom}.
In particular, the notion of \emph{\rum} is defined. In the cases of that paper,
a \rum corresponds to a rigid rotation of $N$ point vortices about some axis.
We then give the equations of motion for point vortices on a rotating sphere \cite{B77} in Section
\ref{motionrot}.
The Hamiltonian system is infinite dimensional due the background vorticity coming from the rotation of the
sphere.

We prove that a \rum formed of latitudinal rings of the non-rotating system persists when
the sphere rotates.
From \cite{LMR00} and \cite{LP02}, we know that the following arrangements of latitudinal rings
are \ria of the non-rotating system:
$C_{Nv}(R,k_pp)$, $D_{Nh}(2R,k_pp)$, $D_{Nd}(R,R',k_pp)$, and $D_{2Nh}(R_e,k_pp)$
($k_p$ is the number of polar vortices).
See Figures \ref{cnvRevpolesch3c}, \ref{fign2}, and \ref{fign1}.
\bfig
\bc
\subfigure[$C_{4v}(R)$]{\label{cnvRch3c}\includegraphics[width=4.6cm,height=5.8cm,angle=0]{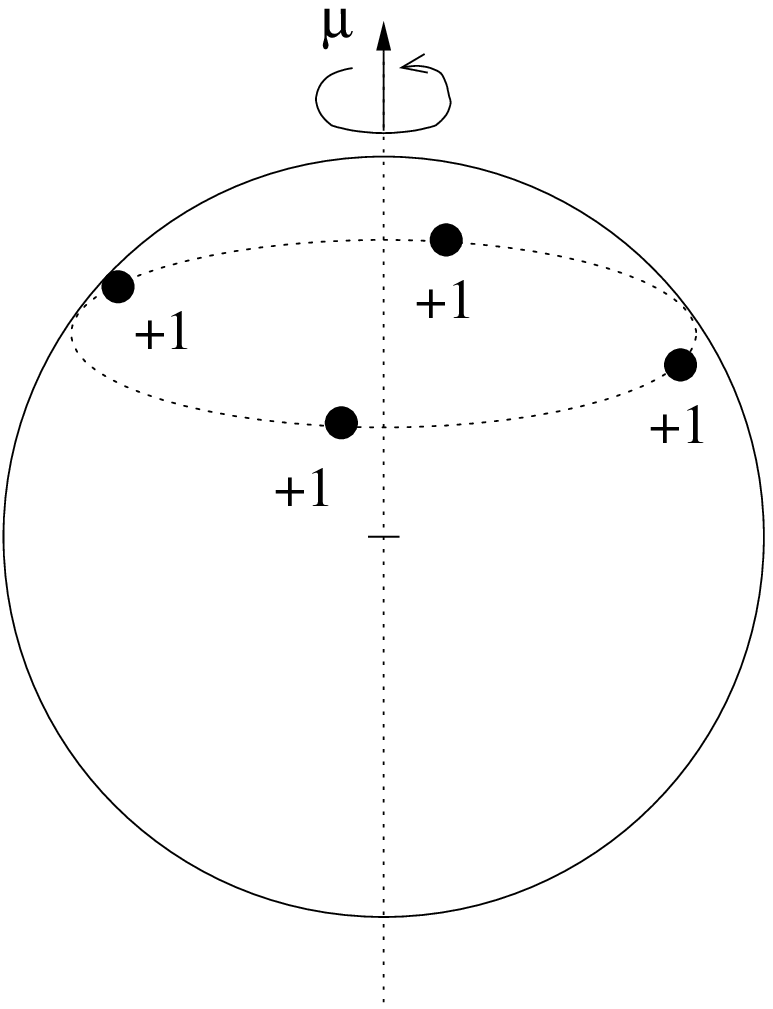}}
\hspace{1.2cm}
\subfigure[$C_{4v}(R,p)$]{\label{cnvRpch3c}\includegraphics[width=4.6cm,height=5.8cm,angle=0]{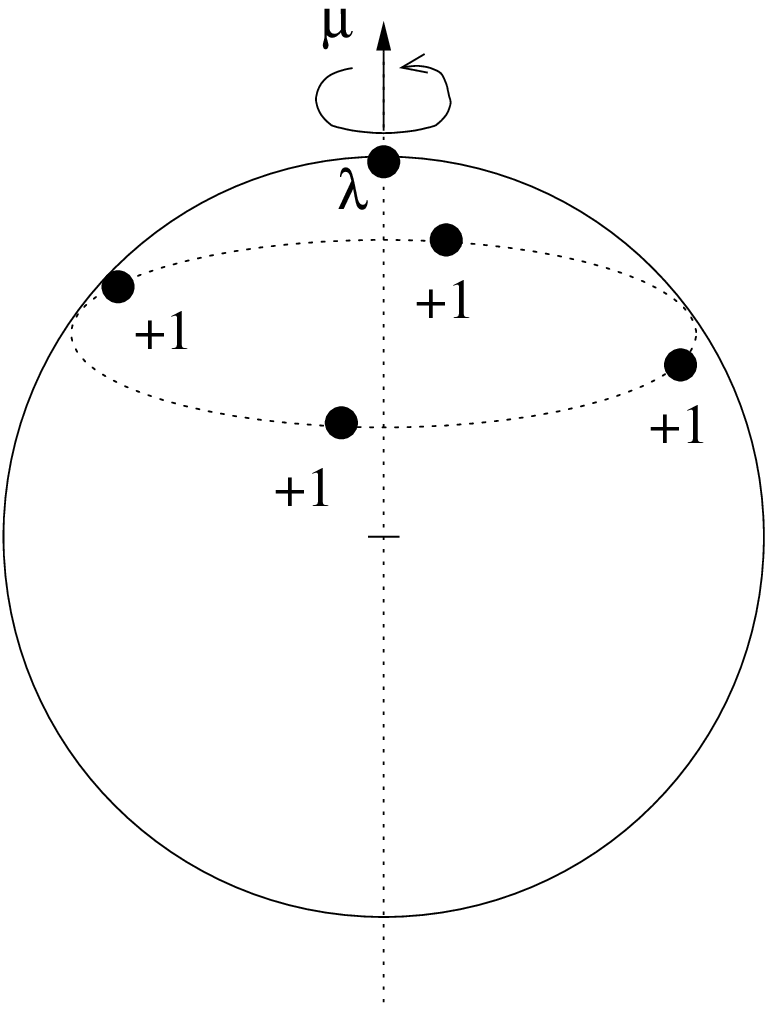}}
\caption{The $C_{nv}(R,k_p p)$ \ria with $k_p=0,1$.}
\label{cnvRevpolesch3c}
\ec
\efig
\bfig
    \bc
\includegraphics{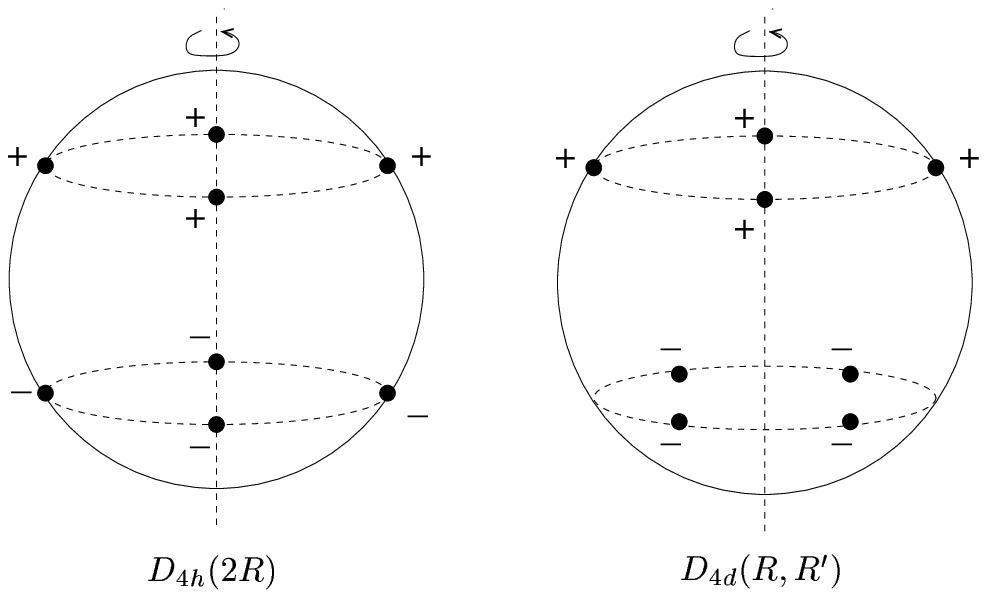}
\caption{Relative equilibria $D_{4h}(2R)$ and $D_{4d}(R,R')$.}
\label{fign2}
     \ec
\efig
\bfig
    \bc
 \includegraphics{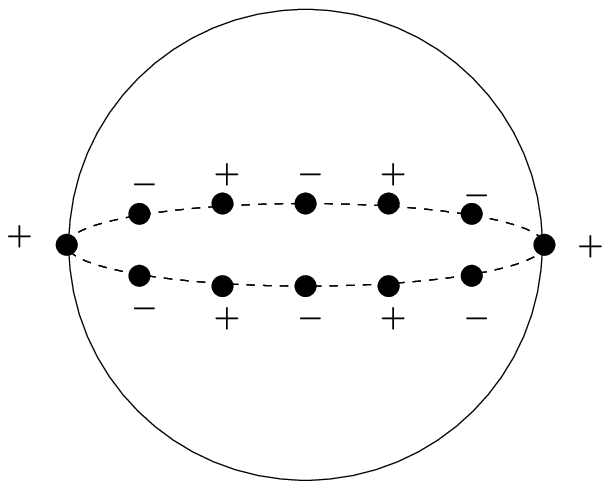}
 \caption{Equatorial ($\pm$)ring $D_{2Nh}(R_e)$.}
 \label{fign1}
     \ec
\efig

In Section \ref{stabrot}, we give the stability of \ria $C_{Nv}(R)$ and $C_{Nv}(R,p)$
for three different approximations or limiting cases:
point vortices on a rotating plane, geostrophic vortices, and point vortices on a non-rotating sphere.
Indeed in that particular cases the system becomes finite dimensional and we can use the
techniques of Section \ref{geom} to obtain both nonlinear and linear stability results.
The stability is determined with a block diagonalization version of the Energy-momentum method
\cite{OR99,LP02} and the Lyapunov stability results are \emph{modulo $SO(2)$}.
In particular, we compute the (nonlinear) stability of a polygon formed of $N$ identical point vortices in
the plane
together with a central vortex of arbitrary vorticity in Appendix B. Our results differ from those
of \cite{CS99} but agree with the linear study of \cite{MS71}.
We also improve some stability results on geostrophic vortices of \cite{MS71} proving that some
linear stable configurations are actually Lyapunov stable.

The paper ends with a discussion on the Southern Hemisphere Circulation and its relationship
with vortices.

\section{Geometric mechanics}
\label{geom}

We recall thereafter some basics notions of geometric mechanics. We refer to \cite{MR94,Or98} for
further details.

Let $G$ be a connected Lie group acting smoothly on a symplectic manifold $(\p,\varpi)$.
Consider an Hamiltonian dynamical system $(\p,\varpi,H)$ with \mom $J:\p\to\g^*$
such that the Hamiltonian vector field $X_H$ and the \mom are $G$-equivariant.
A point $x_e\in\p$ is called a \textit{relative equilibrium} if for all $t$ there exists
$g_t \in G$ such that $x_e(t)=g_t \cdot x_e$, where $x_e(t)$ is the dynamic orbit of $X_H$ with $x_e(0)=x_e$.
In other words, the trajectory is contained in a single group orbit.
A \rum $x_e\in\p$ is a critical point of the \textbf{augmented Hamiltonian}:
$$H_{\xi}(x)=H(x)-\left< J(x)-\mu,\xi\right>$$ for some $\xi\in\g$.
The vector $\xi$ is unique if the action of $G$ is locally free,
and is called the \emph{angular velocity} of $x_e$.

In the case of point vortices on a rotating sphere, we will have $G=SO(2)$.
Hence \ria are rigid rotations with angular velocity $\xi\in\sodeux\simeq\R$ around
the axis of rotation of the sphere.

The Hamiltonian is $G$-invariant, but may have additional symmetries.
Denote by $\hat G$ the group of symmetries of the Hamiltonian.
For example, in the case of $N$ identical point vortices on a non-rotating sphere, one has
$G=SO(3)\times S_N$ and $\hat G=O(3)\times S_N$ \cite{LMR00}.
Let the fixed point set of a subgroup $K$ of $\hat G$ be:
$$
\Fix (K,\p)=\lbrace x\in\mathcal{P} \mid g\cdot x = x, \forall g\in K \rbrace.
$$
The following theorem permits to determine \riab, it is a corollary of the
\textit{Principle of Symmetric Criticality} of Palais \cite{P79}.
\begin{theo}
\label{michou}
Let $K$ be a subgroup of $\hat G$, $x\in\Fix(K,\p)$ and $\mu=\J(x)$.
Assume that $\hat G$ is compact.
If $x$ is an isolated point in $\Fix(K,\p)\cap \J^{-1}(\mu)$, then $x$ is a relative equilibrium.
\end{theo}
Note that this result depends only on the symmetries, the phase space and the momentum map,
and \emph{not} depends on the form of the Hamiltonian.
A relative equilibrium obtained with that result is said to be a
\textit{large symmetry relative equilibrium} since its isotropy subgroup must be large.

To compute the stability of \riab, we use the \textbf{Energy momentum method}:
let $x_e$ be a relative  equilibrium, $\mu=\J(x_e)\in\g^*$ and $\xi$ be its angular velocity.
The method consists first to determine the symplectic slice
$$\G= T_{x_e}(G_{\mu}\cdot x_e)^\bot \cap \Ker D\J(x_e)$$
where
$$
G_{\mu}=\lbrace g\in G\mid \Coad_g\cdot \mu=\mu \rbrace.
$$
The second step consists to examine the definitness of $d^2H_{\xi}|_{\G}(x_e)$,
and apply the following result \cite{Pa92,OR99} which holds in particular if $G$ is compact:

\noindent\textit{If $d^2H_{\xi}|_{\G}(x_e)$ is definite, then $x_e$ is Lyapunov stable modulo $G_{\mu}$.}

In Section \ref{stabplane}, we will consider vortices in the plane (point vortices and geostrophic
vortices), the symmetry group $G$ is $SE(2)$ or $SO(2)$ depending on whether the plane is rotating.
However we will forget translational symmetries since $SE(2)$ is not compact,
hence we take $G=SO(2)$ and the previous result holds.
Moreover we have $G_\mu=SO(2)$ for all $\mu\in\sodeux^*\simeq\R$.
In Section \ref{stabnonrot}, the sphere is non-rotating, thus the symmetry group is $G=SO(3)$
\cite{LP02,LMR}. We have $SO(3)_\mu=SO(2)$ for $\mu\neq0$, and $SO(3)_{\mu=0}=SO(3)$.
Since we will consider only \ria with a non-zero momentum $\mu$ for that section,
\emph{Lyapunov stable} will mean \emph{Lyapunov stable modulo $SO(2)$} throughout that paper.

The linear stability is investigated calculating the eigenvalues of the linearization in the
symplectic slice, that is of $L_\G=\Omega_{\G}^{\flat -1}d^2H_{\xi}|_{\G}(x_e)$
where $\Omega_{\G}^{\flat}$ is the matrix of $\varpi|_{\G}$.

The symplectic slice $\G$ is a $G_{x_e}$-invariant subspace.
Hence we can perform an $G_{x_e}$-isotypic decomposition of $\G$, this permits to block diagonalize
the matrices $d^2H_{\xi}|_{\G}(x_e)$ and $L_\G$, their eigenvalues are then easier to compute and
we can conclude about both Lyapunov and linear stability.
A basis of the symplectic slice in which these matrices block diagonalize is called a
\emph{symmetry adapted basis}. The symmetry adapted bases do not depend on the particular form of
the system, they depend only on the symmetries of the system.

The different steps of the method are widely detailed in \cite{LP02} which is a study on
point vortices on a non-rotating sphere.

A \rum $x_e$ is said to be \emph{elliptic} if it is spectrally stable with $d^2H_{\xi}|_{\G}(x_e)$ not
definite.
An elliptic \rum may be Lyapunov stable, but this can not be proved via the Energy-momentum method
(but KAM theory may work). Note also that an elliptic \rum becomes linearly unstable when some
dissipation is added to the system \cite{DR02}.

\section{Equations of motion and relative equilibria}
\label{motionrot}

In this section we consider $N$ point vortices on a unit sphere rotating with a constant angular
velocity $\fer$ around the axis $(Oz)$. Hence the $SO(3)$ symmetry of the non-rotating system
breaks to $SO(2)$. The Coriolis force induces a continuous vorticity $\om_\fer$ on the sphere:
$\om=\om_0+\om_\fer$ where $\om_0=\sum_i \lambda_i \delta(x-x_i)$, $\lambda_i$ is the vorticity of
the vortex $x_i$, and $\om_\fer(t=0)=2\fer\cos\tet$. The continuous vorticity is not uniform and
thus interacts \emph{a priori} with the singular vorticity i.e. the vortices. This interaction
makes $\om_\fer$ a function of time. For example, without any vortices
$\om=\om_\fer=2\fer\cos\tet$ is a steady solution.

The stream function $\psi$ satisfies $\Delta\psi=\om$. Hence $\psi=G\ast\om$ where $G$ is the
Green function on the sphere $G(x,x^\prime)=1/(4\pi) \ln(1-x\cdot x^\prime)$, and we obtain the
following expression for the stream function: $$
\psi=\frac{1}{4\pi}\sum_{i=1}^{N}\lambda_i\ln(1-x\cdot x_i) + G\ast\om_\fer . $$ The equations
governing the motion of the vortices are therefore \cite{B77}: $$ \dot\theta_i=\sum_{j=1,j\neq i}^{N}
\lambda_j\frac{\sin\theta_j\sin(\phi_i-\phi_j)}{1-x_i\cdot
x_j}+\frac{1}{\sin\theta_i}\frac{\partial G\ast\om_\fer}{\partial\phi_i}(\tet_i,\phi_i) $$ $$
\sin\theta_i\ \dot \phi_i=-\sum_{j=1,j\neq i}^{N} \lambda_j \frac{\sin\theta_i \cos\theta_j
-\sin\theta_j \cos\theta_i \cos(\phi_i-\phi_j)} {1-x_i\cdot x_j}-\frac{\partial
G\ast\om_\fer}{\partial\tet_i}(\tet_i,\phi_i) $$ for all $i=1,\dots,N$. Moreover, the vorticity
satisfies $\partial_t\om+(u\cdot\nabla)\om=0$ (Euler), $u$ being the velocity, and this equation
can be written equivalently with a Poisson bracket: $$ \dot\om=\lbrace \om,\psi \rbrace $$ where
$$ \lbrace f,g \rbrace=\frac{\partial f}{\partial\cos\theta}\frac{\partial
g}{\partial\phi}-\frac{\partial f}{\partial\phi}\frac{\partial g}{\partial\cos\theta} $$ for two
smooth functions $f,g$ on the sphere.
The full dynamical system is therefore: $$ \dot\theta_i=\sum_{j=1,j\neq i}^{N}
\lambda_j\frac{\sin\theta_j\sin(\phi_i-\phi_j)}{1-x_i\cdot
x_j}+\frac{1}{\sin\theta_i}\frac{\partial G\ast\om_\fer}{\partial\phi_i}(\tet_i,\phi_i) $$ $$
\sin\theta_i\ \dot \phi_i=-\sum_{j=1,j\neq i}^{N} \lambda_j \frac{\sin\theta_i \cos\theta_j
-\sin\theta_j \cos\theta_i \cos(\phi_i-\phi_j)} {1-x_i\cdot x_j}-\frac{\partial
G\ast\om_\fer}{\partial\tet_i}(\tet_i,\phi_i) $$ $$\dot\om_\fer=\lbrace \om_\fer,\psi \rbrace$$
since the strengths of the vortices are constant.\\

\rmk The vorticity must satisfy $\int \om\; dS=0$ from Stoke's theorem. It should be noted that
when the sum $\Lambda=\sum\lambda_i$ is non-zero, then the stream function remains as before
taking $\om_0=-\Lambda/(4\pi)+\sum_i \lambda_i \delta(x-x_i)$. Since $\om_0$ satisfies $\int
\om_0\; dS=0$, it follows that $\int \om\; dS=0$ when $\om_\fer=2\fer\cos\tet$.\\

It is easy to verify that the following quantity (total kinetic energy) is conserved and serves as
a Hamiltonian:
\begin{eqnarray*}
H & = & \frac{1}{2}\int_{S^2} \om\cdot\psi\ dS= \frac{1}{2}\int_{S^2} \om\cdot G\ast\om\ dS\\
 & = & \frac{1}{4\pi} \sum_{i<j}\lambda_i\lambda_j\ln(1-x_i\cdot x_j) + \frac{1}{2}\sum_{k=1}^N
\lambda_k G\ast\om_\fer(\tet_k,\phi_k)\\
 &  & +\frac{1}{2}\int_{S^2}\int_{S^2}
\om_\fer(\tet,\phi)\om_\fer(\tet^\prime,\phi^\prime)G(\tet,\phi,\tet^\prime,\phi^\prime)dS
dS^\prime.
\end{eqnarray*}

Let $SO(2)$ be the group of rotations with as axis the axis of rotation of the sphere, and
consider the diagonal action of $SO(2)$ on the product of $N$ spheres. Clearly the continuous
vorticity satisfies $$
\om_\fer(g\cdot\tet,g\cdot\phi,g\cdot\tet_1,g\cdot\phi_1,\dots,g\cdot\tet_N,g\cdot\phi_N)=\om_\fer(\tet,\phi,\tet_1,\phi_1,\dots,\tet_N,\phi_N)
$$ for all $g\in SO(2)$. It follows that $H$ is $SO(2)$-invariant and the dynamical system is
$SO(2)$-equivariant.

Due to the continuous vorticity, the phase space $\p$ becomes infinite dimensional: $$ \mathcal{P}
= \lbrace (x_1,\dots,x_{N}) \in S^2 \times \cdots \times S^2  \mid x_i \neq x_j\ \mbox{if}\ i \neq
j \rbrace \times \mathcal{F}(S^2) $$ where $\mathcal{F}(S^2)$ is the set of smooth functions on
the sphere.

The following vector (momentum vector) is conserved and thus provides three conserved quantities \cite{B77}:
\begin{eqnarray*}
\vec{M} & = & \int_{S^2} \om\cdot\vec{x}\ dS\\
 & = & \sum_{j=1}^N\lambda_j\; x_j + \left(\begin{array}{l} \int_{S^2}\om_\fer(\tet,\phi)\sin\theta\cos\phi\ dS\\
\int_{S^2}\om_\fer(\tet,\phi)\sin\theta\sin\phi\ dS\\ \int_{S^2}\om_\fer(\tet,\phi)\cos\theta\ dS
\end{array}\right).
\end{eqnarray*}
The term $\sum \lambda_j x_j$ corresponds to the moment map for point vortices on a non-rotating
sphere \cite{LMR00}. Indeed in the case of the non-rotating sphere, the dynamical
system is $SO(3)$-equivariant, and this leads to three conserved quantities since $SO(3)$ is of
dimension three. In the case of the rotating sphere the system is only $SO(2)$-equivariant, hence
the symmetries provide only one conserved quantity (which is the $z$-component of $\vec{M}$).
Actually, the three conserved quantities come from a general property of $2D$ incompressible and
inviscid fluid flows on compact and simply connected surfaces which states that the vector $\int
\om\cdot\vec{x}\; dS$ is conserved, no matter the symmetries of the surface we have \cite{B}. That
vector and the momentum map simply coincide in the case of a non-rotating sphere.

We would like to know if the \ria found on the non-rotating sphere persist when the sphere
rotates. We call an \emph{N-ring} a latitudinal regular polygon formed of $N$ identical vortices.
The following theorem show that \ria formed of $N$-rings (together with possibly some polar
vortices)  persist.
\begin{theo}
\label{ria} Let $x_e$ be a \rum of the non-rotating system, with angular velocity $\xi_0$, formed
of $m$ $N_m$-rings ($N_m\geq 2$) together with possibly some polar vortices. Then
$(x_e,\om_\fer(t)=2\fer\cos\tet)$ is a \rum of the rotating system with angular velocity
$\xi=\xi_0+\fer$.
\end{theo}
\begin{coro}
\label{cororia} The \ria
$$C_{Nv}(R,k_p p),\ D_{Nh}(2R,k_pp),\ D_{Nd}(R,R',k_pp),\ D_{2Nh}(R_e,k_pp)$$
persist when the sphere rotates.
\end{coro}
\begin{proof}
Let $(x_e,\xi_0)$ be a \rum of the non-rotating system formed of $k_r$ $N_k$-rings, a polar vortex
being considered as a $1$-ring. Hence
\begin{eqnarray*}
\dot \tet_i & =  & \frac{1}{\sin\theta_i}\frac{\partial
G\ast\om_\fer}{\partial\phi_i}(\tet_i,\phi_i)\\ \dot \phi_i & =  &
\xi_0-\frac{1}{\sin\theta_i}\frac{\partial G\ast\om_\fer}{\partial\tet_i}(\tet_i,\phi_i)
\end{eqnarray*}
for all $i=1,\dots,N$.

Let the continuous vorticity be $\om_\fer(t)=2\fer\cos\tet$ for all time $t$. We have
$\frac{\partial G\ast\om_\fer}{\partial\phi_i}=0$ since $\int G(\tet,,\phi,\tet_i,\phi_i) d\phi$
does not depend on $\phi_i$, and
\begin{eqnarray*}
\frac{\partial G\ast\om_\fer}{\partial\tet_i}&=&\frac{\fer}{4\pi}\int_{\tet=0}^{\pi} \sin 2\tet
\int_{\phi=0}^{2\pi} \frac{ \cos\theta\sin\theta_i -\sin\theta \cos\theta_i
\cos(\phi-\phi_i)}{1-\cos\theta\cos\theta_i-\sin\theta \sin\theta_i \cos(\phi-\phi_i)} d\phi d\tet
\\
 &=&\frac{\fer}{2}\int_{\tet=0}^{\pi} \sin 2\tet \left\lbrack \frac{\cos\theta\sin\theta_i}{|\cos\theta-\cos\theta_i|} +\frac{\cos\theta_i}{\sin\theta_i} \left(1-\frac{1-\cos\theta\cos\theta_i}{|\cos\theta-\cos\theta_i|}\right) \right\rbrack d\tet .
\end{eqnarray*}
And we obtain after some calculus that $\frac{\partial
G\ast\om_\fer}{\partial\tet_i}=-\fer\sin\theta_i$.

Thus we proved that $\dot\tet_i=0$ and $\dot\phi_i=\xi_0+\fer$ for all $i=1,\dots,N$. It remains
then to prove that $\dot\om_\fer=0$. One has $\dot\om_\fer=\lbrace \om_\fer,\psi
\rbrace=2\fer\frac{\partial\psi}{\partial\phi}$. Since $\frac{\partial
G\ast\om_\fer}{\partial\phi}=0$, it follows that $$
\frac{\partial\psi}{\partial\phi}=\sum_{k=1}^{k_r}\sum_{j=1}^{N_k} \lambda_k
\frac{\sin\theta\sin\theta_k\sin(\phi-\phi_{j,k})}{1-\cos\theta\cos\theta_k-\sin\theta
\sin\theta_k \cos(\phi-\phi_{j,k})} $$ where $\lambda_k$ and $\theta_k$ are respectively the
vorticity and the co-latitude of the ring $k$, and $\phi_{j,k}=\epsilon_k+2\pi (j-1)/N$. It can be
shown that this sum vanishes (see Appendix A), hence $\dot\om_\fer=0$ and
$(x_e,\om_\fer(t)=2\fer\cos\tet)$ is a \rum of the rotating system with angular velocity
$\xi=\xi_0+\fer$.
\end{proof}

\rmk In order to prove existence of \riab, we may think to use the \textit{Principle of Symmetric
Criticality} of Section \ref{geom}. The Principle holds since $SO(2)$ is compact,
however the fixed point sets $\Fix (C_N,\p)$ are infinite dimensional, it is therefore
quite a task to find \ria with that method.

\section{Stability of \ria on the rotating sphere}
\label{stabrot}

In this section, we unfortunately do not compute the stability of the \ria determined in the
previous section. Indeed the Hamiltonian system is here infinite dimensional and the method
described in Section \ref{geom} work only for finite dimensional Hamiltonian systems.
Hence we will give in this section the stability of different approximations or limiting cases of
the \pg rotating sphere problem\pd\ which lead to a finite dimensional Hamiltonian system.

\subsection{Vortices on a rotating plane}
\label{stabplane}

We consider point vortices on a rotating plane in order to approximate the $C_{Nv}(R)$ and
$C_{Nv}(R,p)$ \ria on a rotating sphere when the ring is close to the pole. This approximation can
be done in two different manners: the \pg classical\pd\ point vortices on a rotating plane, and the
\pg geostrophic\pd\ vortices \cite{S43}. The equivalent of the $C_{Nv}(R,p)$ arrangement on the plane is
the arrangement $C_{N}(R,p)$: a regular polygon formed of $N$ vortices of strength $1$ together
with a central vortex of strength $\la$. When there is no central vortex, the configuration is
denoted $C_{N}(R)$ and corresponds to the $C_{Nv}(R)$ arrangement on the sphere (see Figure
\ref{cnRevpoles}). We assume that $N\geq 3$ since the cases $N=2$ are degenerate due the
colinearity of the arrangements.\\
\bfig \bc
\subfigure[$C_{4}(R)$]{\label{cnR}\includegraphics[width=4.8cm,height=5cm,angle=0]{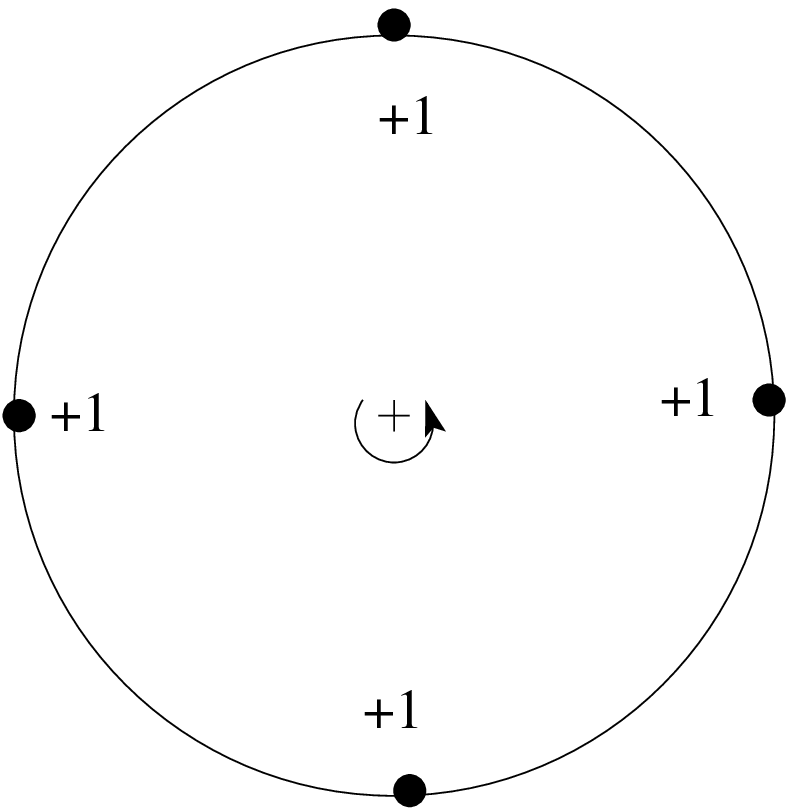}}
\hspace{1.2cm}
\subfigure[$C_{4}(R,p)$]{\label{cnRp}\includegraphics[width=4.8cm,height=5cm,angle=0]{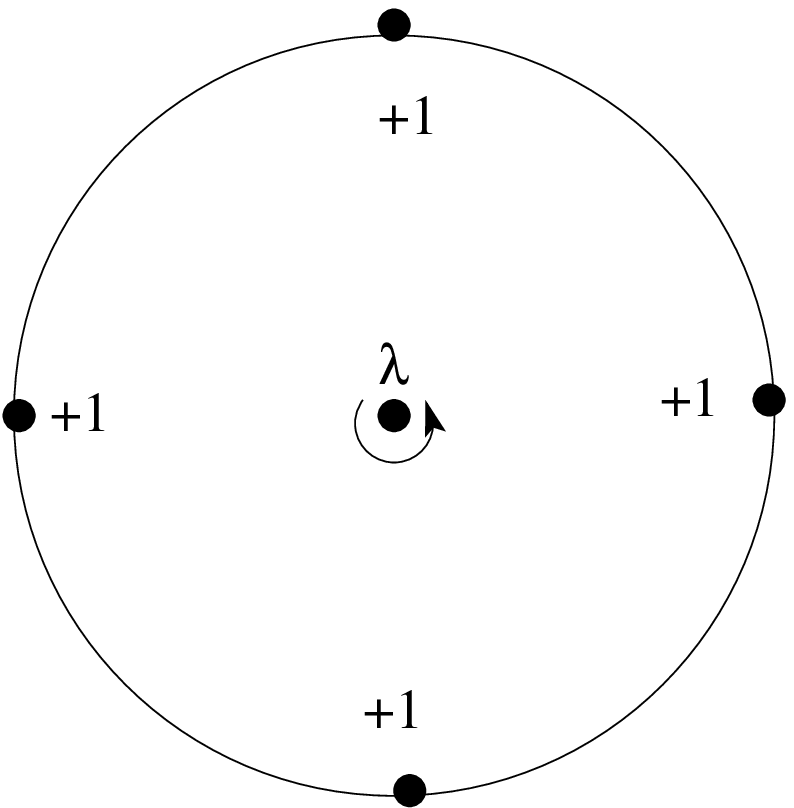}}
\caption{The $C_{N}(R,k_p p)$ planar \riab.} \label{cnRevpoles} \ec \efig

\noindent\textbf{Point vortices on a rotating plane.}\\ Consider a system of point vortices on a
plane rotating with a constant angular velocity $\fer$ around its normal axis. Contrary to the
rotating sphere, the continuous vorticity induced by the rotation is here uniform, thus the
continuous vorticity does not interact with the point vortices. The dynamics of point vortices on
the rotating plane is therefore similar to that for the non-rotating plane as the following
proposition shows.

Denote a \rum $x_e=(z_1,\dots,z_N)\in\C^N$ with angular velocity $\xi$ by $(x_e,\xi)$ and set the
origin of the plane to be the centre of the rotation. Recall that $\Lambda$ is the sum of the
vorticities, that is $\Lambda=\sum_i\lambda_i$.
\begin{theo}
\label{planarpersist} If $(x_e,\xi)$ is a \rum of the non-rotating plane such that $\sum\lambda_j
z_j=0$ and $\Lambda\neq 0$, then $(x_e,\xi+\fer)$ is a \rum of the rotating plane. Moreover
$(x_e,\xi+\fer)$ satisfies the stability properties of $(x_e,\xi)$ (as elements of $\p/SO(2)$).
\end{theo}
\begin{proof}
The continuous vorticity induced by the rotation of the plane is $\omega_\fer(x,y)=2\fer$, thus
the distribution of vorticity for $N$ point vortices $z_j=\rho_j\exp(i\phi_j)$, $j=1,\dots,N$ is
given by $\omega=2\fer+\sum_j \la_j\ \delta (z-z_j)$. It follows that the Hamiltonian is
$$H=\frac{\fer}{2}\sum_j  \la_j\rho_j^2 + H_0,$$
where $H_0$ is the Hamiltonian for the non-rotating plane (see Appendix B).

Since $\la_j\dot\rho_j=-\frac{1}{\rho_j}\frac{\partial H}{\partial\phi_j}$ and
$\la_j\dot\phi_j=\frac{1}{\rho_j}\frac{\partial H}{\partial\rho_j}$, the dynamical system is given
by $$ \left( \begin{array}{c} \dot\rho_j \\ \dot\phi_j
\end{array}\right)
=\left( \begin{array}{c} 0 \\ \fer
\end{array}\right)+
\left( \begin{array}{c} -\frac{1}{\lambda_j\rho_j}\frac{\partial H_0}{\partial\phi_j}\\
\frac{1}{\lambda_j\rho_j}\frac{\partial H_0}{\partial\rho_j}
\end{array}\right).
$$ Since $\sum\lambda_j z_j=0$ and $\Lambda\neq 0$, one can show that the \rum of the non-rotating
plane $(x_e,\xi)$ satisfies $$ \frac{\partial H_0}{\partial\phi_j}=0,\
\frac{1}{\lambda_j\rho_j}\frac{\partial H_0}{\partial\rho_j}=\xi $$ for all $j$ ($\Lambda\neq 0$
is required to insure the existence of the \pg barycentre\pd\ $\sum_j\lambda_j z_j/\Lambda$).
Hence $(\tilde x_e,\tilde\xi)$ with $\tilde x_e=x_e$ and $\tilde\xi=\fer+\xi$ is a \rum of the
rotating plane.

If a vortex $z_j$ of the \rum $x_e$ is at the origin (a central vortex), this discussion is not
valid due to the degeneracy of polar coordinates at the origin, we then use cartesian coordinates
and find the same conclusions.

We are interested in stability in $\p/SO(2)$ (stability modulo $SO(2)$) where $SO(2)$ is the group
of rotations around the origin. The \pg rotating\pd\ dynamical system is $SO(2)$-equivariant. The
momentum map coming from the $SO(2)$ symmetry is $J=\sum_j
\la_j\rho_j^2/2\in\R\simeq\sodeux^\star$ (see Appendix B), thus $H=H_0+\fer J$. Hence $$
d^2H_{\tilde\xi}(\tilde x_e)=d^2H(x_e)-(\fer+\xi)d^2J(x_e)=d^2H_0(x_e)-\xi
d^2J(x_e)=d^2{H_{0}}_{\xi}(x_e). $$ Moreover, the symplectic slice is the same for both the
\emph{rotating} and the \emph{non-rotating} cases, since the stability is investigated in
$\p/SO(2)$ (and not in $\p/SE(2)$ as we could do for the non-rotating case). The stability
statement follows.
\end{proof}

\rmk The previous result is actually due to the fact that the Hamiltonian $H_0$ is perturbed by a
function $J^\fer(x)=\< J(x),\fer\m$ where $J$ is a momentum map of a symmetry group $G$ of $H_0$,
and $\fer\in\g$. The Hamiltonian $J^\fer$ is a \emph{collective} Hamiltonian (that is a function
of $J$), see \cite{GS84} for a general study of collective Hamiltonians.\\

Consequently, all \ria of the non-rotating plane formed of regular polygons persist to be \ria
when the plane rotates, provided that the polygons are cocentric with centre of the rotation. In
particular, the $C_{N}(R)$ and $C_{N}(R,p)$ \ria persist (the centre of the ring being the centre
of the rotation). For \ria formed of several regular polygons, we refer to \cite{LR96} in which a
stability study is also given for some particular configurations.

We can now give stability results for the $C_{N}(R)$ and $C_{N}(R,p)$ \ria thanks to the stability
study of $C_{N}(R)$ \ria \cite{CS99} on the non-rotating plane and a theorem proved in Appendix B.

\begin{prop}
The $C_{N}(R)$ \ria on a rotating plane are Lyapunov stable if $N\leq 7$ and linearly unstable
otherwise.
\end{prop}

Let $\varepsilon_N=+1$ if $N$ is even, and $\varepsilon_N=0$ if $N$ is odd.
\begin{theo}
\label{cabral} The $C_{N}(R,p)$ \ria with $N\geq4$ on a plane (rotating or not) are Lyapunov
stable if $\max(0,(N^2-8N+7+\varepsilon_N)/16) < \la < (N-1)^2/4$, elliptic if
$(N^2-8N+7+\varepsilon_N)/16 < \la < 0$, and linearly unstable if the previous conditions are not
satisfied. The $C_{3}(R,p)$ \ria are Lyapunov stable if $0<\la<1$, elliptic if $\la<0$, and
linearly unstable if $\la>1$.
\end{theo}

We have $(N^2-8N+7+\varepsilon_N)/16<0$ if and only if $N\leq6$, hence there exists some linearly
stable configurations with a negative central vorticity only for $N\leq6$.\\

\begin{proof}
The stability results of $C_{N}(R,p)$ \ria on a non-rotating plane are given in Appendix B since
the proof is quite long. These results will hold for a rotating plane by Theorem
\ref{planarpersist}.
\end{proof}

\rmk The Lyapunov stability of $C_{7}(R)$ (the Thomson heptagon) is not proved with an energetic
method. Indeed the Hessian is not definite at this \rumb, one needs to go to the fourth order to
determine the Lyapunov stability \cite{M78,CS99}.\\

\noindent\textbf{Geostrophic vortices.}\\ In his paper \cite{S43}, Stewart considered a single
layer of fluid of constant density on a rotating disc in order to model the atmosphere of the
Earth. Let $\fer$ be the constant angular velocity of the disc, $z$ be the complex coordinate on
the rotating disc, $h_0$ be the uniform depth of the fluid at rest and  $g$ be the acceleration
due to gravity. By means of the geostrophic wind equations, he found that the stream function
corresponding to a vortex at $z_0$ with vorticity $\lambda$ is $$ \Psi=\la K_0(\kappa |z-z_0|) $$
where $\kappa=2\fer/\sqrt{gh_0}$ and $K_0$ is the Bessel function of the second kind. Such a
vortex is called a \emph{geostrophic} (or \emph{Bessel}) vortex. The Hamiltonian for a system of
$N$ geostrophic vortices $z_1,\dots,z_N$ with vorticities $\la_1,\dots,\la_N$ is $$ H=\sum_{i<j}
\la_i\la_j K_0(\kappa |z_i-z_j|). $$ Note that a geostrophic vortex is almost a planar
$\kappa^{-1}$-Euler vortex, the parameter $\kappa^{-1}$ is the  horizontal length scale of the
geostrophic vortex. In the case of the Earth $\kappa^{-1}$ is about $2000$ kms. If we rescale the
radius of the Earth to be $1$ as the spheres previously considered, $\kappa$ is about $3$.
Note also that the limiting case $\kappa=0$ corresponds to logarithmic vortices, that is to
\emph{classical} point vortices.

Since the Hamiltonian for the geostrophic vortices has the same symmetries as the Hamiltonian for
the planar point vortices, we follow Section \ref{geom} to determine large symmetry \ria and
study their stability: it is easy to show that $C_{N}(R)$ and $C_{N}(R,p)$ configurations are
large symmetry \ria for the system of planar point vortices, and to find the symmetry adapted
basis. The symmetry adapted basis for $C_{N}(R,p)$ is given in Appendix B, while the one for
$C_{N}(R)$ is easily deduced from the former. These bases are similar to those for the $C_{Nv}(R)$
and $C_{Nv}(R,p)$ configurations of point vortices on a sphere \cite{LMR}.

We set the value of the vorticity of the ring to be $1$. The stability results are given in
Proposition \ref{stabR} for $C_{N}(R)$ geostrophic \ria and in Figure \ref{stabRpfig} for
$C_{N}(R,p)$ geostrophic \riab. The linear stability results were already found by Morikawa and
Swenson \cite{MS71} by means of a numerical computation of the eigenvalues of the linearization.
Here we compute numerically the eigenvalues of the diagonal blocks of the Hessian%
\footnote{But the Hessian is obtained analytically.}
 and the linearization, and thus conclude about both formal (hence Lyapunov) and linear stability.
The results of this paper agree with their results except for the right stability frontiers where
they differ slightly. This is probably due to a lack of accuracy. The advantage of the method
described here is that we know that both the Hessian and the linearization block diagonalize,
hence it is not necessary to compute the nil blocks (gain of time), and we know that the
components in it are all zero and this contributes to improve the accuracy of the eigenvalues
(indeed, if we compute numerically the components of the nil blocks, we find for example
$10^{-13}$ instead of zero, and this leads finally to a lack of accuracy in the computations of
the eigenvalues).

\begin{prop}
\label{stabR} The stability of the geostrophic $C_{N}(R)$ \ria depends on $N$ and $\kappa$ as
follows:
\begin{itemize}
\item $N\geq 7$: linearly unstable for all $\kappa>0$;
\item $N=6$: Lyapunov stable for $\kappa\in [0,1.28]$, linearly unstable otherwise;
\item $N=5$: Lyapunov stable for $\kappa\in [0,3.75]$, elliptic otherwise;
\item $N=3,4$: Lyapunov stable for all $\kappa\geq 0$.
\end{itemize}
\end{prop}
\bfig \bc
\includegraphics[width=5.6cm,angle=0]{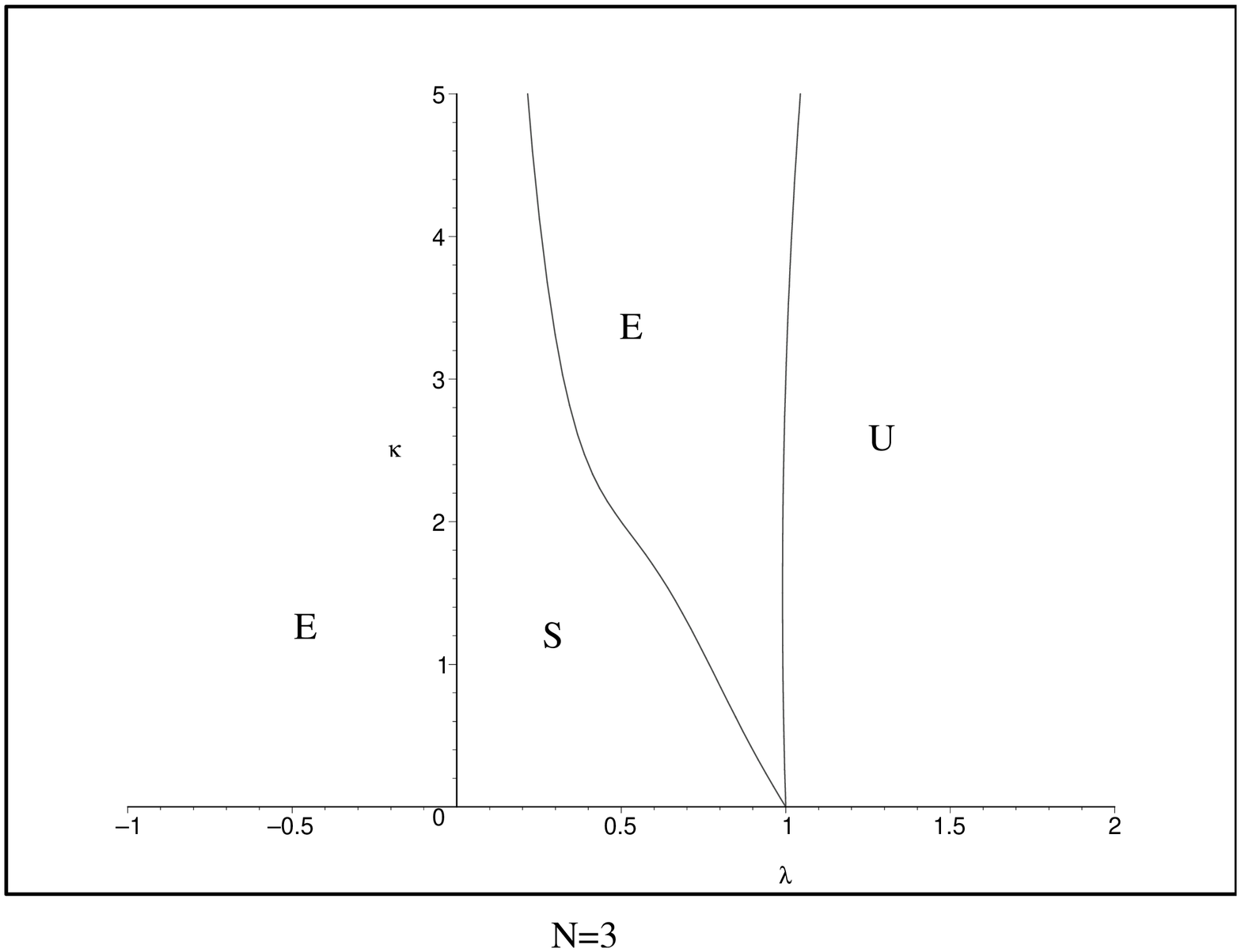}
\includegraphics[width=5.6cm,angle=0]{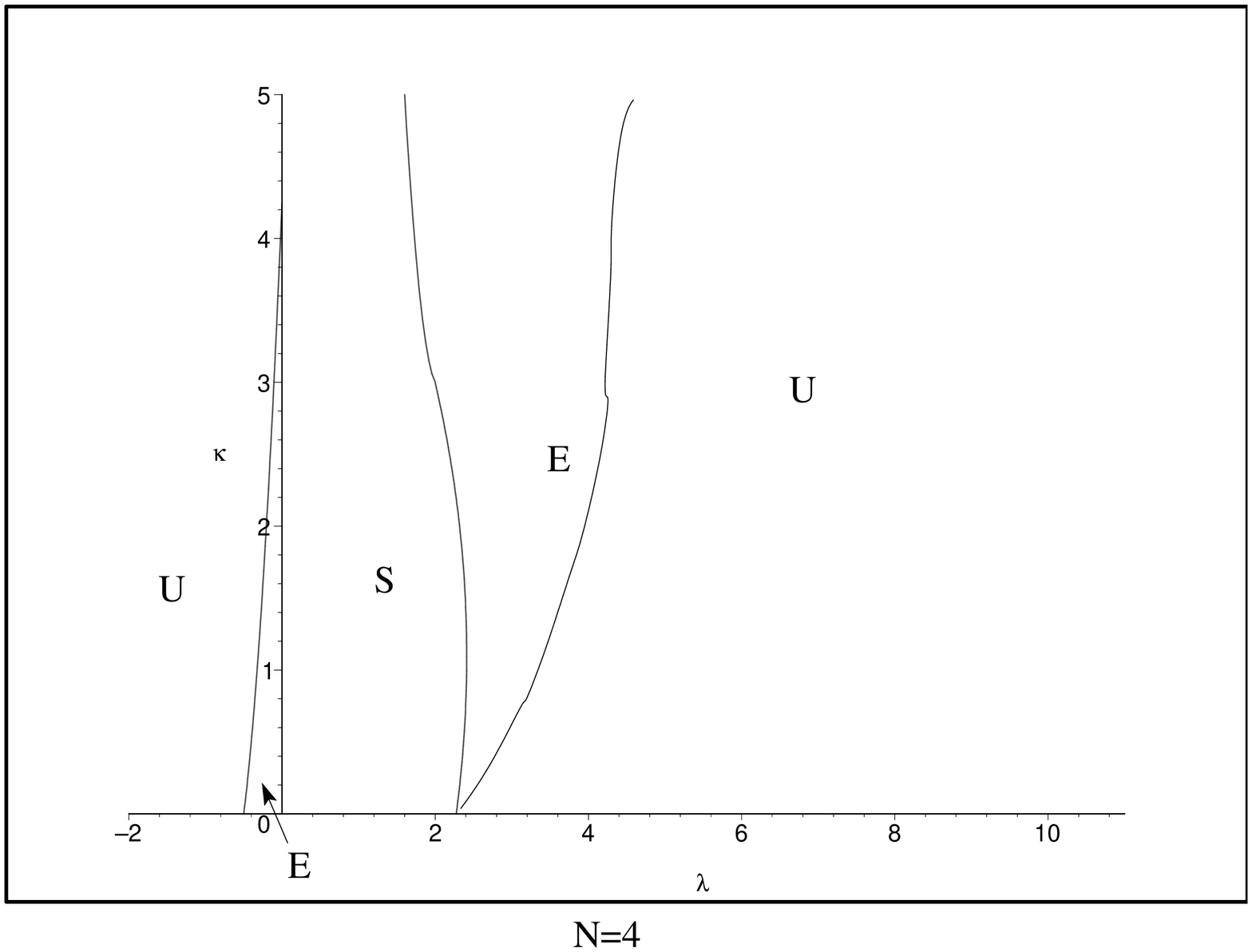}

\vspace{0.5cm}

\includegraphics[width=5.6cm,angle=0]{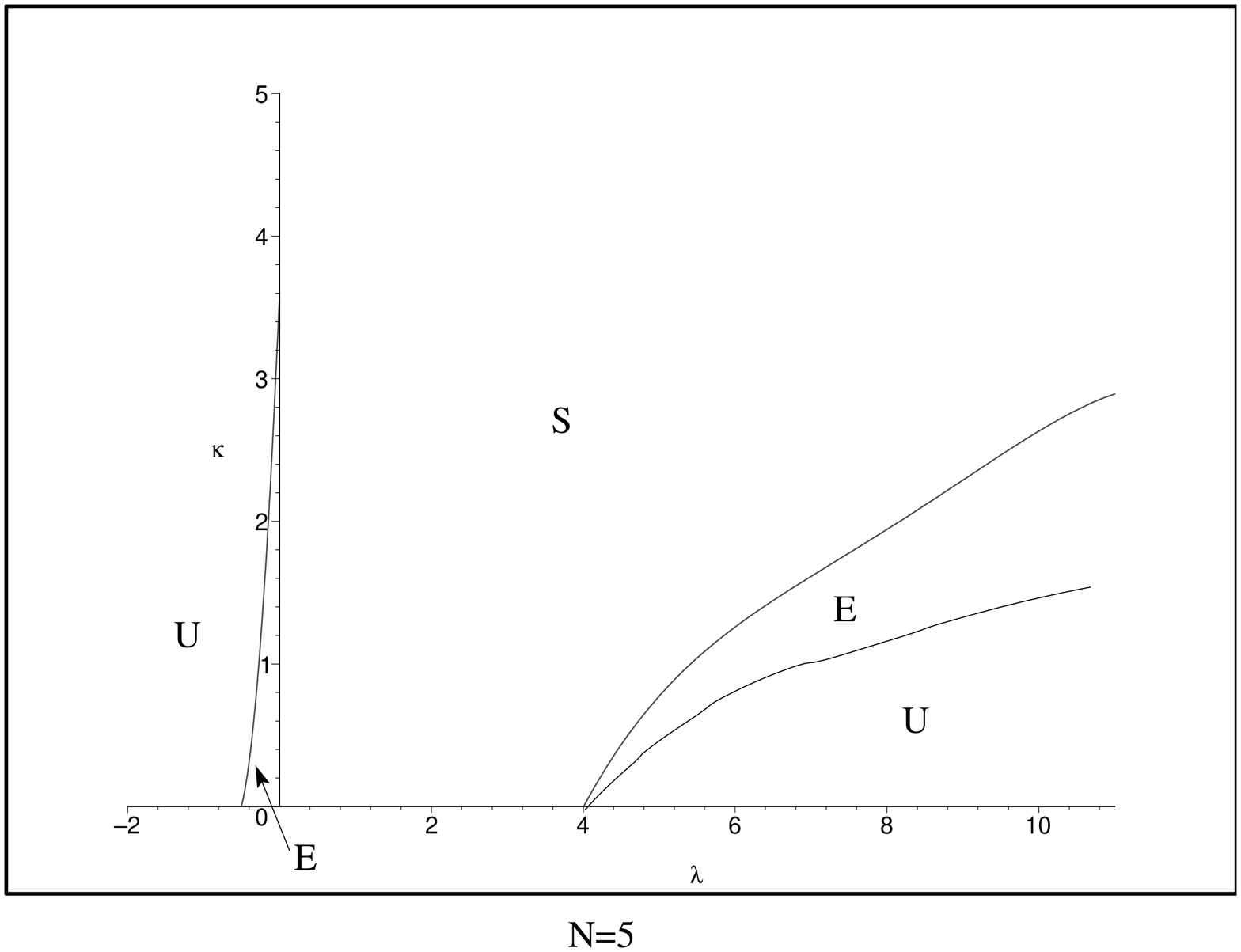}
\includegraphics[width=5.6cm,angle=0]{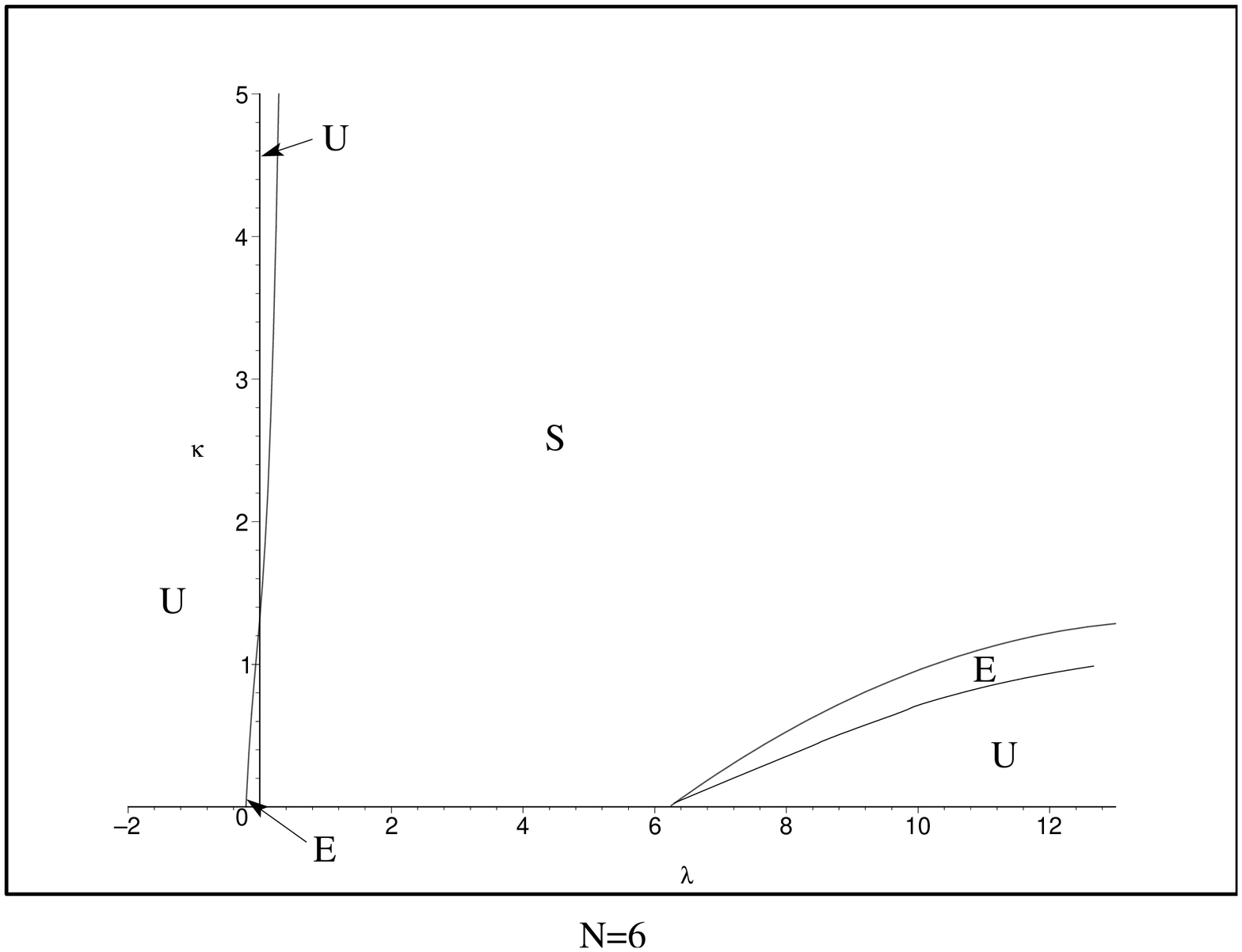}

\vspace{0.5cm}

\includegraphics[width=5.6cm,angle=0]{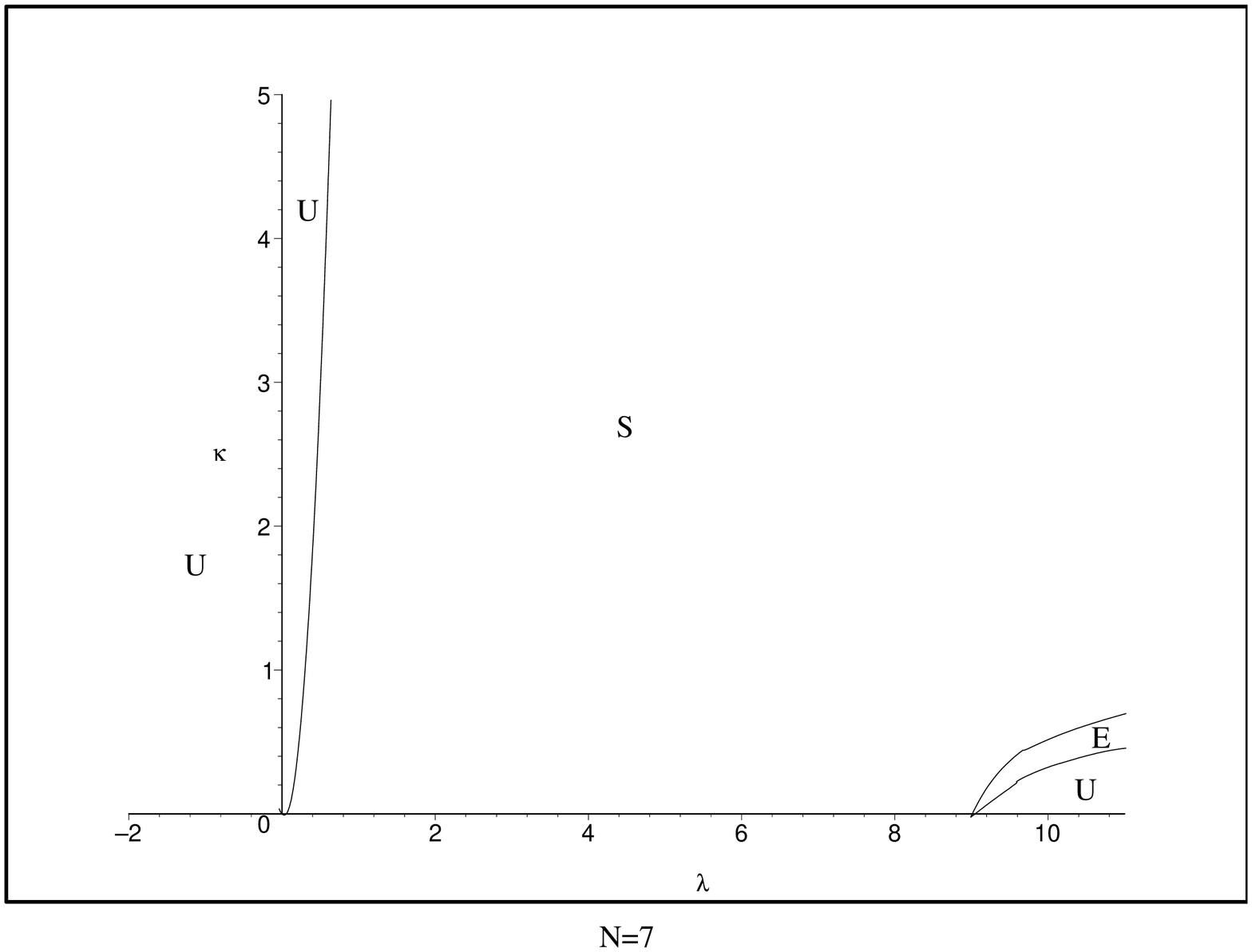}
\includegraphics[width=5.6cm,angle=0]{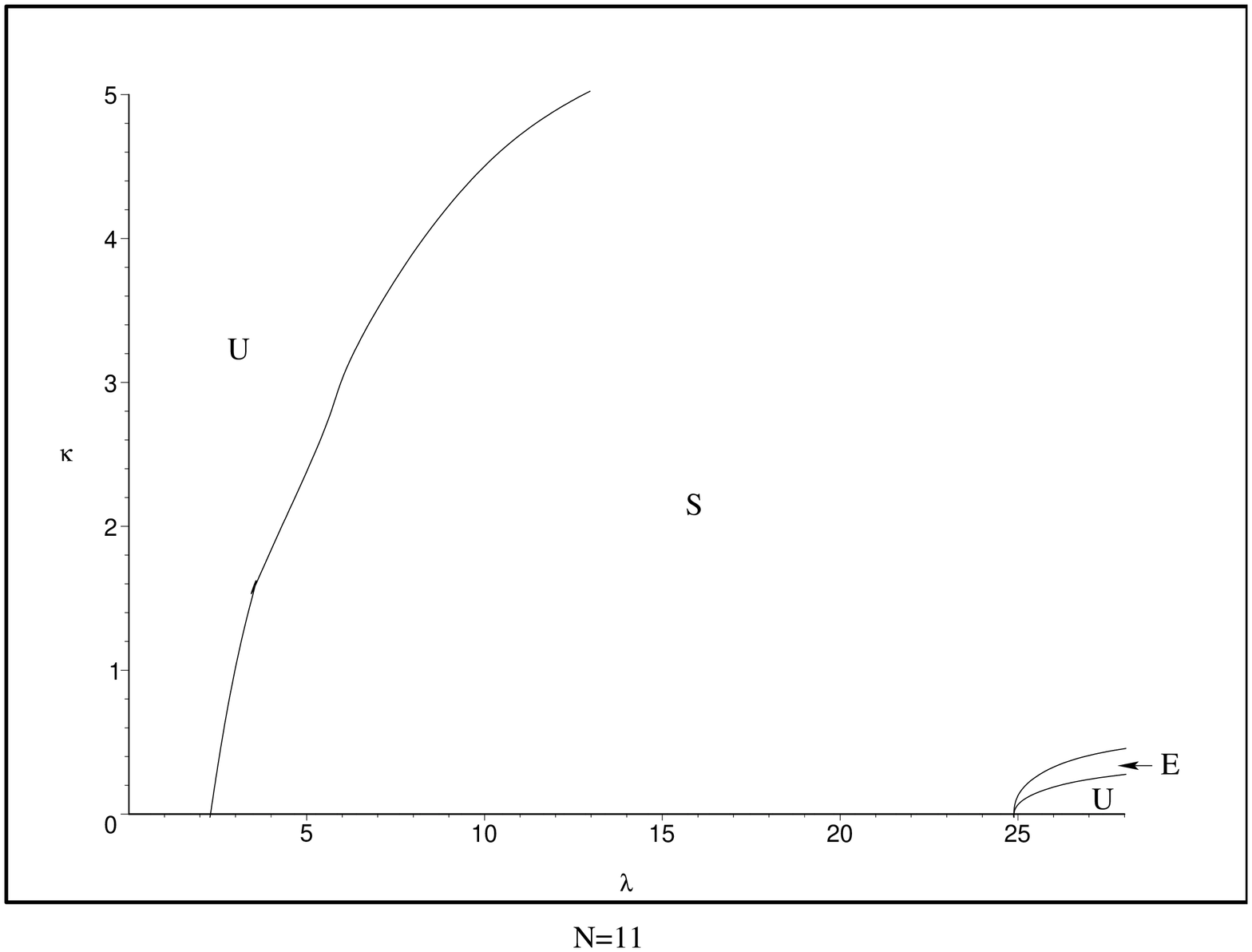}
     \caption{Stability results for geostrophic $C_{N}(R,p)$ \riab.
              A \rum can be formally stable (S), linearly unstable (U), elliptic (E).
              The diagrams for $N\geq 8$ are similar to that for $N=11$.}
     \label{stabRpfig}
\ec \efig

Instability on the left part of the figures is due to the ring, while instability in right part is
due to the central vortex. For $N=3$, the linear stability frontier of the central vortex (the
right one) is almost a vertical line for $\kappa\in[0,5]$, that is $\kappa$ in that range has
almost no influence on the linear stability of the central vortex.
For $N=4,\dots,7$, the stability frontier of the ring (the left one) has the same property. On the
other hand, the stability frontiers of the central vortex (the two on the right) for $N>4$ depend
on $\kappa$. Note that the notion of \pg dependence\pd\ is linked to the scale of the figure for
$N=5$, that is $(\kappa,\lambda)\in[0,5]\times [-2,10]$.

For $N\leq7$, almost every \rum stable for $\kappa=0$ is stable for $\kappa\in [0,5]$ . Thus if a
point vortices $C_{N}(R,p)$ \rum is stable, then so does the corresponding geostrophic \riab. The
point vortices system is therefore a good model to determine the stability (not instability) for
$N\leq7$. This becomes less and less true as $N$ increases. Indeed the greater is $N$, the more
the stability depends on $\kappa$.

\subsection{Stability results for $\Omega=0$}
\label{stabnonrot}

In this section, we give the stability results for the $C_{Nv}(R)$ and $C_{Nv}(R,p)$ \ria on the
\emph{non-rotating} sphere.

The following theorem was proved first in term of linear stability in \cite{PD93}, and recently
in term of Lyapunov stability in \cite{LMR} and \cite{BC01} by different ways.

\begin{theo}
  The stability of the ring of $N$ identical vortices depends on $N$
  and the latitude as follows:
  \begin{description}
  \item[N=2] is Lyapunov stable at all latitudes;
  \item[N=3] is Lyapunov stable at all latitudes;
  \item[N=4] is Lyapunov stable if $\cos^2\theta_0 > 1/3$, and linearly unstable if the inequality is reversed;
  \item[N=5] is Lyapunov stable if $\cos^2\theta_0 > 1/2$, and linearly unstable if the inequality is reversed;
  \item[N=6] is Lyapunov stable if $\cos^2\theta_0 > 4/5$, and linearly unstable if the inequality is reversed;
  \item[N$>$6] is always linearly unstable,
  \end{description}
  where $\theta_0$ is the co-latitude of the ring.
\end{theo}

\begin{figure}[htbp]
    \begin{center}
\includegraphics[width=5.6cm,angle=0]{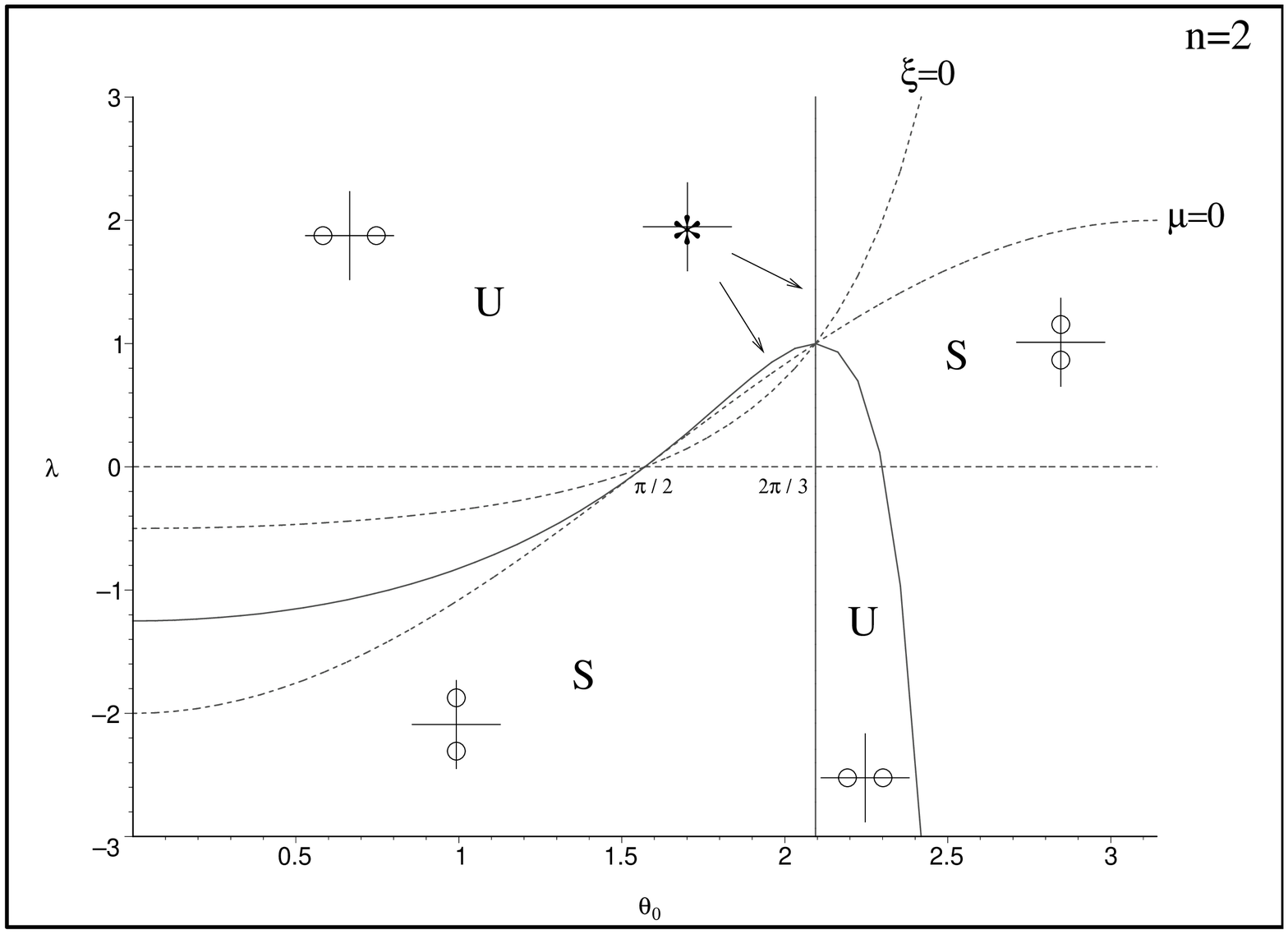}
\includegraphics[width=5.6cm,angle=0]{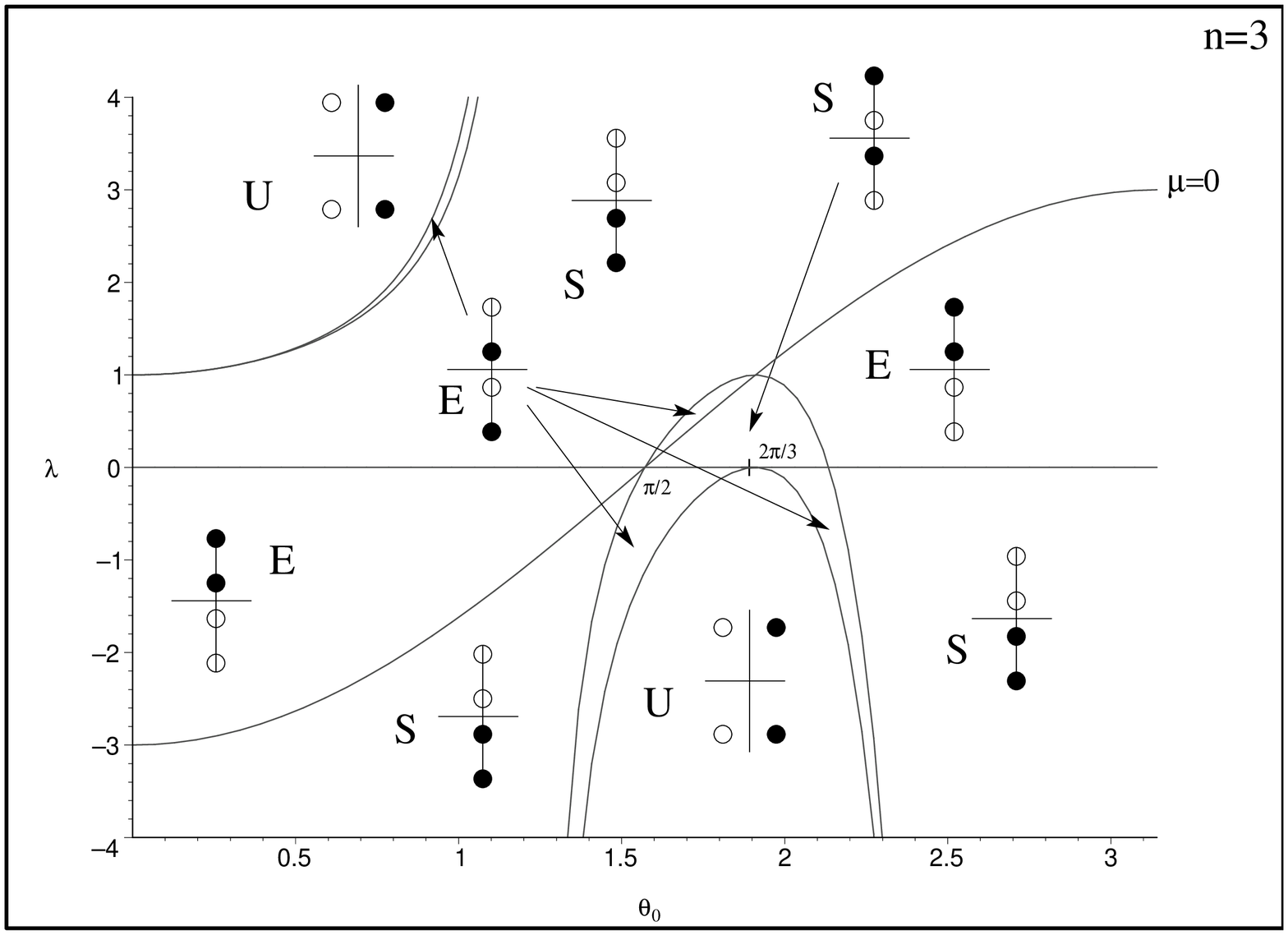}
\caption{Bifurcation diagrams for  $C_{2v}(R,p)$ and $C_{3v}(R,p)$ relative equilibria.
     The circles represent the eigenvalues of the mode $1$.
     A \rum can be Lyapunov stable (S), linearly unstable (U), elliptic (E).}
\label{CRp1}
     \end{center}
 \end{figure}

\begin{figure}[htbp]
    \begin{center}
\includegraphics[width=5.6cm,angle=0]{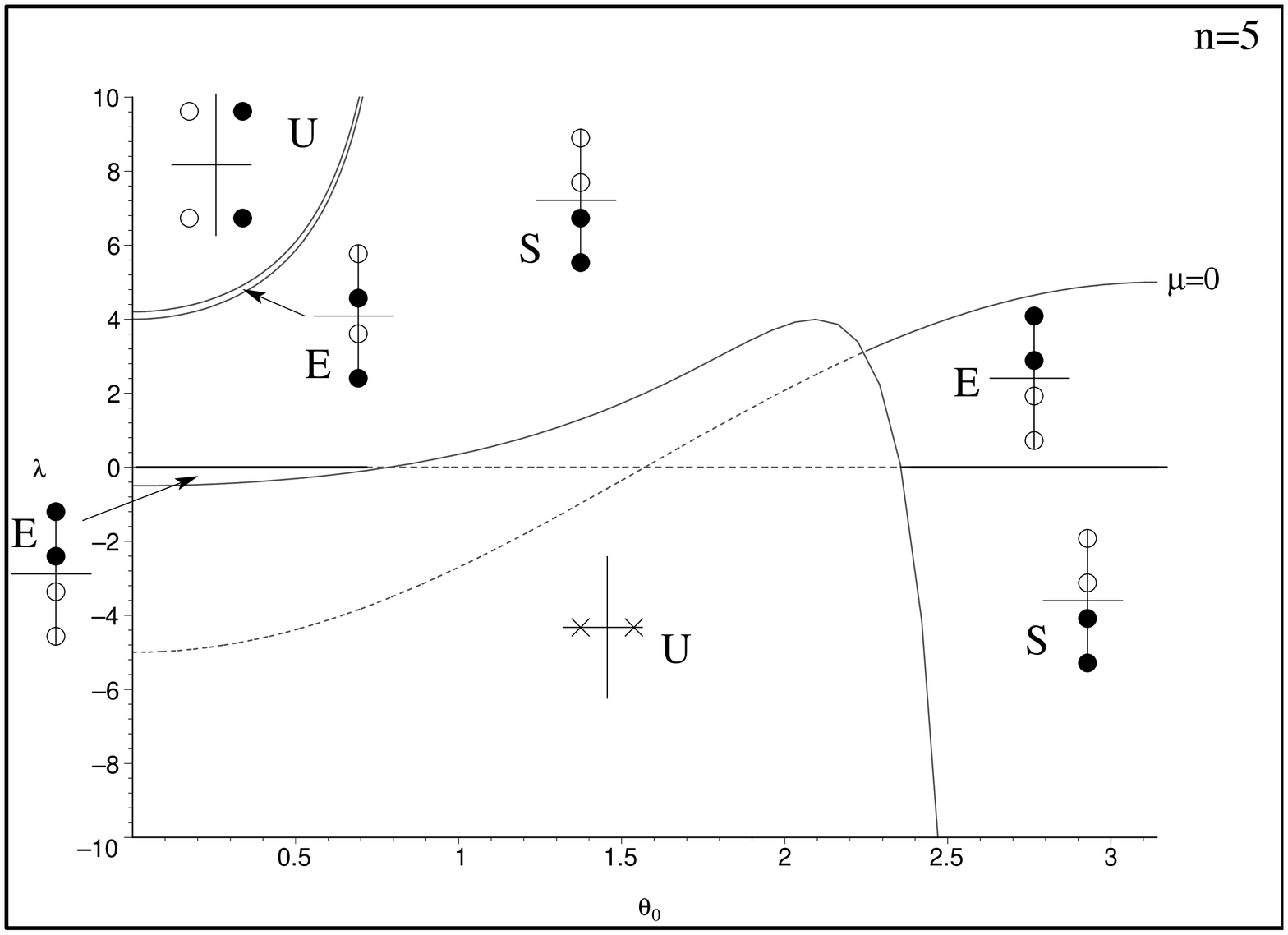}
\includegraphics[width=5.6cm,angle=0]{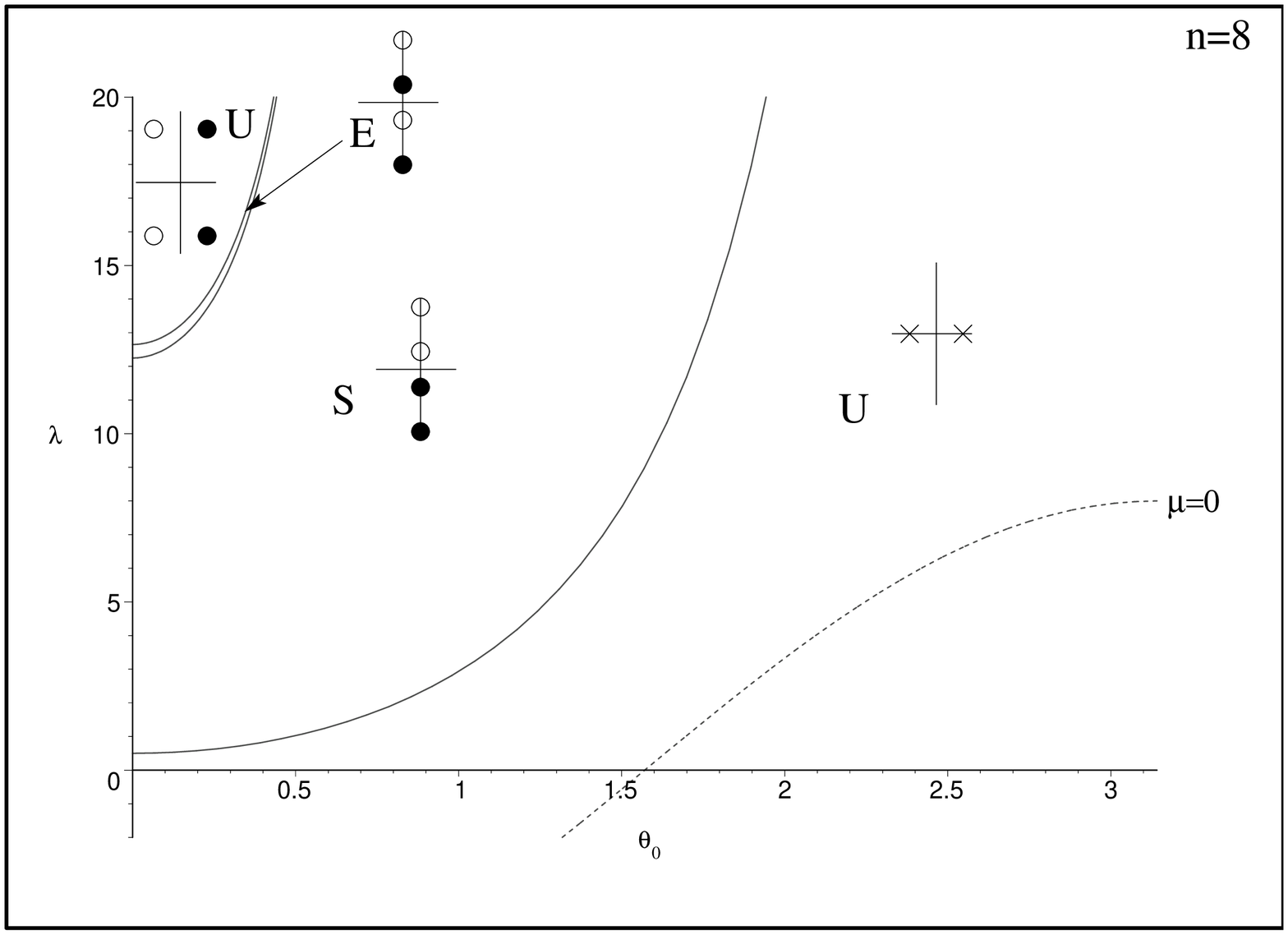}
     \caption{Bifurcation diagrams for  $C_{Nv}(R,p)$ relative equilibria.
     The bifurcation diagrams for $N\geq7$ are similar to that for
     $N=8$, while those for $N=4,6$ are similar to that for $N=5$.
     The circles represent the eigenvalues of the mode $1$, while the cross represent those of the mode $[N/2]$.
     A \rum can be Lyapunov stable (S), linearly unstable (U), elliptic (E).}
\label{CRp2}
     \end{center}
 \end{figure}

The stability of $C_{Nv}(R,p)$ \ria is given in the next theorem and illustrated by Figures \ref{CRp1}
and \ref{CRp2}. We assume that the
polar vortex is at the North pole without loss of generality, and that the momentum of the \rum is
non-zero. Let $\theta_0$ be the co-latitude of the ring.

\begin{theo}[Laurent-Polz, Montaldi, Roberts \cite{LMR}]\ \par\noindent
A $C_{2v}(R,p)$ \rum is Lyapunov stable if
$$(1+2\cos\theta_0)[(1+\cos\theta_0)^2\la+\cos\theta_0(2+3\cos\theta_0)]<0. $$ and linearly
unstable if the inequality is reversed.\\ A $C_{3v}(R,p)$ \rum is Lyapunov stable if $a\lambda
(\lambda+n\cos\theta_0)(\lambda-\lambda_1) < 0$, and spectrally unstable if and only if $$
8a\lambda >(n\sin^2\theta_0+4(n-1)\cos\theta_0)^2. $$ A $C_{Nv}(R,p)$ relative equilibrium with
$N\geq 4$

\noindent (i) is spectrally unstable if and only if $$ \lambda < \lambda_0\ \ \mbox{or}\ \
8a\lambda > (N\sin^2\theta_0+4(N-1)\cos\theta_0)^2, $$ \noindent (ii) is Lyapunov stable if $$
\lambda > \lambda_0\ \ \mbox{and}\ \ a\lambda (\lambda+N\cos\theta_0)(\lambda-\lambda_1) < 0, $$
where $$
\begin{array}{lll}
a&=&(N\cos\theta_0-N+2)(1+\cos\theta_0)^2\\ \lambda_1&=& (N-1)\cos\theta_0\;
(N\sin^2\theta_0+2(N-1)\cos\theta_0)/a\\ \lambda_0&=&(c_N -
(N-1)(1+\cos^2\theta_0))/(1+\cos\theta_0)^2\\ c_N&=&\cases{N^2/4 & if $N$ is even, \cr
              (N^2-1)/4 & if $N$ is odd.}
\end{array}
$$
\end{theo}

\bigskip

\ \par\noindent Let $x_e$ be a $C_{Nv}(R)$ or $C_{Nv}(R,p)$ \rum on a sphere rotating with angular
velocity $\fer$. By the theory of perturbations, there exists a neighborhood $U$ of zero,
$U\subset\R$, such that for all $\fer\in U$, $x_e$ has the stability of the corresponding \rum on
the non-rotating sphere. Hence the stability results of the two previous theorem persist in a
neighborhood of zero for $\fer$.


\subsection{Conclusions}

From the results of that section, one can reasonably think that for a small angular velocity
$\fer$ the \rum $C_{Nv}(R)$ with a small co-latitude $\theta_0$ (that is the ring is close to the
pole) is Lyapunov stable for $N\leq5$, while it is linearly unstable for $N\geq8$. The cases $N=6$
and $N=7$ are less clear since the stability for $N=6$ depends on the parameter $\kappa$ in
Proposition \ref{stabR}, and the case $N=7$ is stable or unstable depending on whether we are on
the plane or the sphere.

About the \rum $C_{Nv}(R,p)$ with a negative vorticity for the polar vortex and a small
co-latitude $\theta_0$ (that is the ring and the polar vortex vorticities of have opposite signs),
one can reasonably think from the previous results that for a small angular velocity $\fer$ the
\rum is:
\begin{itemize}
  \item elliptic for $N=3$,
  \item elliptic for a weak polar vorticity and $N=4,5,6$,
  \item linearly unstable otherwise.
\end{itemize}

\section{The Southern Hemisphere Circulation}
\label{southern}

It is well known from sailors that around Antarctica there are all year long strong winds coming
from the West, the \emph{westerlies}.
The atmospheric circulation around Antarctica is indeed essentially \emph{zonal} during all seasons of
the year \cite{RL53,S68}, the global circulation is a steady rotation of low pressure systems
around the North-South axis. This circulation is commonly called \emph{The Southern Hemisphere
Circulation} by meteorologists. With the words of this paper, we can say that the atmospheric
circulation is close to the \rum $C_{Nv}(R)$ formed of $N$ low pressure systems. Since there
exists a polar vortex (a high pressure system) at the South pole, the relevant \rum is merely
$C_{Nv}(R,p)$. The southern hemisphere circulation varies with seasons \cite{L65,L67} but also
intraseasonally \cite{RL53,MF91}: the westerlies are stronger in summer, and their latitude may
fluctuate within a season, while the zonal character is always maintained. Polar high pressure
systems, which have their origin in or near Antarctica, at times attain sufficient strength to
break through the zonal flow to reach and reinforce a sub-tropical high pressure system (the
anticyclonic system of Saint-Helen for example).

Can point vortices explain the Southern Hemisphere Circulation? Indeed for sixty years some
meteorologists or physicists study point vortices in order to understand motions of low/high
pressure systems in the Earth's atmosphere \cite{S43,MS71,F75,DP98}. In \cite{S43} the linear
stability of the geostrophic $C_3(R)$ is proved and according to the author, that result explains
the presence of the three sub-tropical major high pressure systems on each hemisphere. For example
in the southern hemisphere, we have the anticyclonic systems of Saint-Helen, of the Indian Ocean,
and of the South Pacific, which persist. In the northern hemisphere, the high pressure systems of
the Acores is well known from europeans. In \cite{F75}, a linear stability study of $C_{3v}(R)$ in
the $\beta$-plane approximation is used to conclude that a ring of major high pressure systems is
stable at sub-tropical latitude, while for low pressure systems it does not. However, a ring of
low pressure systems is stable near the poles. The problem of those two papers is that
sub-tropical latitudes are partially occupied by continents, and they do not take in account the
continent-atmosphere interaction.

In \cite{DP98}, the authors relate their work on vortices to atmospheric blocking events. The
advantage of the Southern Hemisphere Circulation is that it occurs around Antarctica, which is a
continent with a disc shape, thus the $SO(2)$ symmetry is not broken. Moreover, the change
continent-ocean on the Antarctica coast must influence the atmospheric circulation, but --- again
--- this does not break the $SO(2)$ symmetry. Note also that this region (40°-90° South) is the
only one large region of the Earth which has an $SO(2)$ symmetry.

From the previous section, one can conjecture that a ring of $N$ low pressure systems is Lyapunov
stable for $N=3,4,5$; and a ring of $N=3,4,5,6$ low pressure systems together with a weak polar
high pressure system is elliptic, while the whole system becomes unstable if the polar high gets
stronger. This agrees with the observations of the beginning of this section. But first, we did
not take into account the sub-tropical highs comparing the southern hemisphere circulation to
$C_{Nv}(R,p)$ \riab. Indeed the sub-tropical highs may balance the influence of the polar high,
and therefore reinforce the stability of the ring of low pressure systems. Second, we did not take
into account the external heating by the Sun which adds dissipation to the system, and we know
that dissipation generically induces instability for elliptic \riab. However, this instability may
be compensated by other phenomena such as the presence of sub-tropical highs. With the external
heating by the Sun, the vorticity is no longer a conserved quantity, and one prefers the
\emph{potential vorticity} which is almost conserved (see \cite{DL93}). The reader may have a look
on the forecast for the southern atmospheric circulation on the web at the following address:
\texttt{http://www.ecmwf.int/products/forecasts/d/charts/deterministic/world/}.

\vspace{1.5cm}

\noindent {\large\textbf{Acknowledgements}}\\
\\
This work on point vortices on a rotating sphere was suggested by James Montaldi and Pascal
Chossat.
I would like to gratefully acknowledge and thank James Montaldi for very helpful comments and suggestions,
and Guy Plaut for interesting discussions on meteorology.

\newpage\noindent{\Large\textbf{Appendix A}}

\bigskip

\noindent Let $B_k$ be the sum $$ B_k=\sum_{j=1}^{N_k}
\frac{\sin(\phi-\phi_{j,k})}{1-r\cos(\phi-\phi_{j,k})} $$ where $\phi_{j,k}=\epsilon_k+2\pi
(j-1)/N$, $r=\sin\theta \sin\theta_k/(1-\cos\theta\cos\theta_k)$ and $\phi$ is fixed. We have $$
B_k=\sum_{j=1}^{N_k} \sum_{m=0}^{\infty} r^m \sin(\phi-\phi_{j,k})\cos^m(\phi-\phi_{j,k}) $$ for
all $|r|< max \lbrace 1/|\cos(\phi-\phi_{j,k})| , j=1,\dots,N_k \rbrace$. Since there exist
$\al_0,\dots,\al_m$ such that $\cos^m(\phi-\phi_{j,k})=\sum_{p=0}^{m}
\al_p\cos(p(\phi-\phi_{j,k}))$, we have $B_k=\sum_{m=0}^{\infty} r^m C_m$ where $$
C_m=\sum_{p=0}^{m} \al_p \sum_{j=1}^{N_k}\sin(\phi-\phi_{j,k})\cos(p(\phi-\phi_{j,k})). $$ It
follows from $\sum_{j=1}^{N_k}\sin(q(\phi-\phi_{j,k}))=0$ that $C_m=0$ for all $m$, hence $B_k=0$
for all $k$.

\vspace{2cm}\noindent{\Large\textbf{Appendix B}}

\bigskip

\noindent We give in this appendix the proof of the stability of the $C_{N}(R,p)$ \ria on a
non-rotating plane, that is the proof of the following theorem:\\
\emph{Let $\varepsilon_N=+1$ if $N$ is even, and $\varepsilon_N=0$ if $N$ is odd.\\ The
$C_{N}(R,p)$ \ria with $N\geq4$ on a non-rotating plane are Lyapunov stable if
$\max(0,(N^2-8N+7+\varepsilon_N)/16) < \la < (N-1)^2/4$, elliptic if $(N^2-8N+7+\varepsilon_N)/16
< \la < 0$, and linearly unstable if the previous conditions are not satisfied. The $C_{3}(R,p)$
\ria are Lyapunov stable if $0<\la<1$, elliptic if $\la<0$, and linearly unstable if $\la>1$.}\\
The proof is quite long though we omitted lots of details. We first recall some well-known facts
about the dynamics of planar point vortices.

Let $x_e$ be a $C_{N}(R,p)$ \rum and $\lambda$ be the vorticity of the central vortex. The
equations of motion $N+1$ planar point vortices are \cite{Ar82}:
$$\dot{\overline{z_j}}=\frac{1}{2\pi i}\sum_{k=1,k\neq j}^{N+1} \frac{\lambda_k}{z_j-z_k}$$ where
$z_j$ is a complex number representing the position of the $j$-vortex (we identified the plane
with $\C$), and  $\lambda_j$ is the vorticity of the $j$-vortex. The Hamiltonian for this system
is $$ H=-\frac{1}{4\pi}\sum_{i<j} \lambda_i \lambda_j \ln\vert z_i-z_j \vert^2 $$ and the symmetry
group of the vector field $X_H$ is $SE(2)$ (which is not compact). In the particular case of $N$
identical vortices together with an additional vortex, the Hamiltonian is $E(2)\times
S_N$-invariant. Let $z_j=\rho_j\exp (i\phi_j),j=1,\dots,N$ and $z_{N+1}=x+iy$, the dynamical
system in $(\rho_j,\phi_j,x,y)$ variables is:
$$\begin{array}{lll}\la_j\dot\rho_j&=&-\frac{1}{\rho_j}\frac{\partial H}{\partial\phi_j}\\
\la_j\dot\phi_j&=&\frac{1}{\rho_j}\frac{\partial H}{\partial\rho_j}\end{array} ,\
\begin{array}{lll}\la\dot x&=&-\frac{\partial H}{\partial y}\\
\la\dot y&=&\frac{\partial H}{\partial x}.\end{array}$$ After identifying $SE(2)$ with $\C\rtimes
S^1$ and so $\se^*$ with $\C\times\R$, the \mom of the system is $$\J=\left( i\sum_{j=1}^{N+1}
\lambda_j z_j , \frac{1}{2}\sum_{j=1}^{N+1} \lambda_j \vert z_j \vert^2 \right).$$ In a frame such
that $z_{N+1}=0$, a $C_N(R,p)$ configuration has coordinates $\rho_j=R, \phi_j=2\pi j/N, j=\1N,
x=0, y=0$, and is a \rum with angular velocity $\xi=(N-1+2\lambda)/4R^2$ (where $H$ was normalized
by $4\pi$). It is straightforward to see that $R^2H_\xi$ does not depend on $R$, hence so does the
stability, and we can set $R=1$.

The recent results of Patrick, Roberts and Wulff \cite{PRW01} generalize the Energy-Momentum
method
 to Hamiltonian systems with a non-compact symmetry group. However, we will apply
here the classical Energy-Momentum method with the compact subgoup $SO(2)$, and thus forget the
translational symmetries. The \mom coming from the rotational symmetries is
$\J=\la(x^2+y^2)/2+\sum_{j}\lambda_j \rho_j^2/2$. Since the coadgoint action of $SO(2)$ is
trivial, one has $SO(2)_\mu=SO(2)$. The symplectic slice at a $C_{N}(R,p)$ \rum is therefore
$$\NN=Span\left\lbrace\sum_{i=1}^N\delta\rho_i,\sum_{i=1}^N\delta\phi_i\right\rbrace^\bot.$$

Let $G=O(2)\times S_N$. The linear map associated to $d^2H_{\xi}|_{\G}(x_e)$ (see Introduction) is
equivariant under both symplectic and anti-symplectic elements of $G_{x_e}$, while $L_\G$  is
equivariant under only the symplectic ones (see \cite{LP02} and \cite{LMR} for a detailed
account on the sphere). We then perform an isotypic decomposition and find that the symmetry
adapted basis of $\NN$ is $$\left( \alpha^{(1)}_{\phi},\beta^{(1)}_{\tet},\widetilde{\delta y},
\beta^{(1)}_{\phi},\alpha^{(1)}_{\tet},\widetilde{\delta x},\left\lbrace
\alpha^{(\ell)}_{\theta},\alpha^{(\ell)}_{\phi},\beta^{(\ell)}_{\theta},\beta^{(\ell)}_{\phi} \mid
2\leq \ell\leq [n/2] \right\rbrace  \right) $$ where $$
\begin{array}{rcl}
\alpha^{(\ell)}_{\theta} +i\beta^{(\ell)}_{\theta} &=& \sum_{s=1}^N\exp(2i\pi\ell
s/N)\delta\rho_{s}\\[4pt] \alpha^{(\ell)}_{\phi} +i\beta^{(\ell)}_{\phi} &=&
\sum_{s=1}^N\exp(2i\pi\ell s/N)\delta\phi_{s}\\[4pt] \widetilde{\delta x}&=&
-\frac{N}{\sqrt{2N}}\delta x\\ \widetilde{\delta y}&=&-\frac{N}{\sqrt{2N}}\delta y.
\end{array}$$ In this basis, the matrix $d^2H_{\xi}|_{\G}(x_e)$ (resp. $L_\G$) block diagonalizes in
$2\times2$ blocks with two $3\times3$ blocks (resp. $4\times4$ blocks with one $6\times6$ block,
plus a $2\times2$ block if $N$ is even).

We first study the Lyapunov stability. Some simple calculations show that $$
\begin{array}{lll}
d^2H_\xi(x_e)\cdot (\beta^{(\ell)}_{\theta},\beta^{(\ell)}_{\theta})&=&d^2H_\xi(x_e)\cdot
(\alpha^{(\ell)}_{\theta},\alpha^{(\ell)}_{\theta})\\
d^2H_\xi(x_e)\cdot(\beta^{(\ell)}_{\phi},\beta^{(\ell)}_{\phi}) &=&d^2H_\xi(x_e)\cdot
(\alpha^{(\ell)}_{\phi},\alpha^{(\ell)}_{\phi})\\
d^2H_\xi(x_e)\cdot(\beta^{(\ell)}_{\theta},\beta^{(\ell)}_{\phi})&=&d^2H_\xi(x_e)\cdot
(\alpha^{(\ell)}_{\theta},\alpha^{(\ell)}_{\phi})=0\\
 d^2H_\xi(x_e)\cdot
(\alpha^{(1)}_\phi,\widetilde{\delta y})&=&d^2H_\xi(x_e)\cdot
(\beta^{(1)}_{\tet},\widetilde{\delta y})\ =\ d^2H_\xi(x_e)\cdot
(\beta^{(1)}_{\phi},\widetilde{\delta x})\\
 &=&d^2H_\xi(x_e)\cdot(\alpha^{(1)}_{\tet},\widetilde{\delta x})=N\la\sqrt{N/2}\\
d^2H_\xi(x_e)\cdot (\widetilde{\delta y},\widetilde{\delta y})&=&d^2H_\xi(x_e)\cdot
(\widetilde{\delta x},\widetilde{\delta x})=N\xi\la/2.
\end{array}
$$ Hence it follows from the block diagonalization that $$ d^2H_\xi|_\NN(x_e)=diag(A,A,D) $$ where
$D=diag(\lbrace\lambda^{(\ell)}_\tet,\lambda^{(\ell)}_\phi,\lambda^{(\ell)}_\tet,\lambda^{(\ell)}_\phi\mid
2\leq \ell\leq [n/2]\rbrace )$, $$ A=\left( \begin{array}{ccc} \lambda^{(1)}_\phi&0&a_N\\
0&\lambda^{(1)}_\tet&a_N\\ a_N&a_N&N\xi\la/2
\end{array}\right)
$$ and $$
\begin{array}{lll}
\lambda^{(\ell)}_\tet&=&d^2H_\xi(x_e)\cdot (\alpha^{(\ell)}_{\theta},\alpha^{(\ell)}_{\theta})\\
\lambda^{(\ell)}_\phi&=&d^2H_\xi(x_e)\cdot (\alpha^{(\ell)}_{\phi},\alpha^{(\ell)}_{\phi})\\
a_N&=&N\la\sqrt{N/2}.
\end{array}
$$ Note that $D$ exists only if $N\geq 4$.

Thanks to the formula \cite{H75} $$\sum_{j=1}^{N-1} \frac{\cos(2\pi\ell j/N)}{\sin^2(\pi
j/N)}=\frac{1}{3}(N^2-1)-2\ell(N-\ell)$$ we find after some lengthy computations that
$\lambda^{(\ell)}_\phi=N\ell(N-\ell)/2$ and $$\lambda^{(\ell)}_\tet=\frac{N}{2}\lbrack
-(\ell-1)(N-\ell-1)+N-1+4\lambda\rbrack .$$ The eigenvalues $\lambda^{(\ell)}_\phi$ are all
positive, thus $D$ is definite if $-(\ell-1)(N-\ell-1)+N-1+4\lambda > 0$ for all
$\ell=2,\dots,[N/2]$, that is if $\lambda > (([N/2]-1)(N-[N/2]-1)-N+1)/4$ which corresponds to
$\lambda > (N^2-8N+8)/16$ for $N$ even and to $\lambda > (N^2-8N+7)/16$ for $N$ odd.

The relative equilibrium is therefore Lyapunov stable if $A$ is positive definite, that is if
 the three following subdeterminants are positive:
$$\lambda^{(1)}_\phi ,\ \lambda^{(1)}_\phi\lambda^{(1)}_\tet,\ \det A.$$ Since
$\lambda^{(1)}_\phi=N(N-1)/2$, $\lambda^{(1)}_\tet=N(4\la+N-1)/2$ and $$\det
A=-N^3\la\left(\la-\frac{N-1}{2}\right)\left(\la-\frac{(N-1)^2}{4}\right),$$ it follows that $A$
is positive definite if and only if $0<\la<(N-1)^2/4$.

We proved therefore that $C_{N}(R,p)$ \ria with $N\geq4$ are Lyapunov stable modulo $SO(2)$ if
$\max(0,(N^2-8N+7+\varepsilon_N)/16)<\la<(N-1)^2/4$ where $\varepsilon_N=0$ (resp.
$\varepsilon_N=1$) for $N$ odd (resp. even).

We then study the linear stability of the relative equilibrium. It follows from the block
diagonalization of $L_\NN$ that $L_\NN=diag\left(A_L, \left\lbrace A_\ell \mid 2\leq \ell\leq
[n/2]\right\rbrace \right)$ where $$ A_L=\left( \begin{array}{cccccc} 0 & 0 & 0 & N\xi/2 & a_N/\la
& a_N/\la\\ 0 & 0 & 0 & -a_N & 0 & -\lambda^{(\ell)}_\tet\\ 0 & 0 & 0 & -a_N &
-\lambda^{(\ell)}_\phi & 0\\ -N\xi/2 & -a_N/\la & -a_N/\la & 0 & 0 & 0\\ a_N & 0 &
\lambda^{(\ell)}_\tet & 0 & 0 & 0\\
 a_N & \lambda^{(\ell)}_\phi & 0 & 0 & 0 & 0
\end{array}\right)
$$ and $$A_\ell=\left( \begin{array}{cccc} 0 & 0 & 0 &-\lambda^{(\ell)}_\phi\\ 0 & 0 &
\lambda^{(\ell)}_\tet & 0\\ 0 &-\lambda^{(\ell)}_\phi & 0 & 0\\ \lambda^{(\ell)}_\tet &  0 & 0 & 0
\end{array}\right) \ \mbox{while it is}$$
$$\left(\begin{array}{cc} 0 &-\lambda^{(\ell)}_\phi\\ \lambda^{(\ell)}_\tet & 0
\end{array}\right) \ \mbox{for}\ \ell=N/2\ (N \mbox{even})$$ (the blocks are given up to a
factor).

After some calculations, we obtain that the eigenvalues of $L_\NN$ are $$ \pm\frac{N}{2}\xi
i,\pm\sqrt{\la-\frac{(N-1)^2}{4}},\left\lbrace\pm i
\sqrt{\lambda^{(\ell)}_\tet\lambda^{(\ell)}_\phi}\mid 2\leq \ell\leq [n/2]\right\rbrace .$$ We
have therefore some double eigenvalues, but the Jordan forms of the blocks are semi-simple. It
follows from the discussion on Lyapunov stability that the eigenvalues $\pm
i\sqrt{\lambda^{(\ell)}_\tet\lambda^{(\ell)}_\phi}$ are all purely imaginary if and only if
$\lambda \geq (N^2-8N+7+\varepsilon_N)/16$. The Theorem follows for $N\geq4$.

For $N=3$, the Lyapunov (resp. linear) stability is governed only by $A$ (resp. $A_L$). Hence
$C_{3}(R,p)$ is Lyapunov stable if $0<\la<(N-1)^2/4$ and linearly unstable if and only if
$\la>(N-1)^2/4$.


\begin{thebibliography}{99}
\bibitem[Ar82] {Ar82} H. Aref, Point vortex motion with a center of symmetry.
                      \emph{Phys. Fluids} {\bf 25} (1982), 2183-2187.
\bibitem[B] {B} G. Batchelor, An Introduction to Fluid Dynamics.
                     Cambridge University Press (1967).
\bibitem[B77]{B77} V. Bogomolov, Dynamics of vorticity at a sphere. \emph{Fluid Dyn.} {\bf 6} (1977), 863-870.
\bibitem[B85]{B85} V. Bogomolov, On the motion of a vortex on a rotating sphere.
                   \emph{Izvestiya} {\bf 21} (1985), 298-302.
\bibitem[BC01]{BC01} S. Boatto and H. Cabral, Non-linear stability of \ria of vortices on a non-rotating sphere.
                   \emph{Preprint Bureau des Longitudes, Paris} (2001).
\bibitem[CS99]{CS99} H. Cabral and D. Schmidt, Stability of \ria in the problem of $N+1$ vortices.
                     \emph{SIAM J. Math. Anal.} {\bf 31} (1999), 231-250.
\bibitem[DL93] {DL93} D. Dritschel and B. Legras,  Modeling oceanic and atmospheric vortices.
                       \emph{Physics Today}  (March 1993), 44-51.
\bibitem[DP98] {DP98} M. Dibattista and L. Polvani,  Barotropic vortex pairs on a rotating sphere.
                       \emph{J. Fluid Mech.} {\bf 358} (1998), 107-133.
\bibitem[DR02] {DR02} G. Derks and T. Ratiu,  Unstable manifolds of relative equilibria in Hamiltonian
                      systems with dissipation. \emph{Nonlinearity} {\bf 15} (2002), 531-549.
\bibitem[F75] {F75} S. Friedlander,  Interaction of vortices in a fluid on a surface of a rotating sphere.
                       \emph{Tellus XXVII} (1975), 1.
\bibitem[GS84]{GS84} V. Guillemin and S. Sternberg, Symplectic Techniques in Physics.
                      Cambridge University Press (1984).
\bibitem[H75] {H75} E. Hansen, A Table of Series and Products.
                      Prentice-Hall (1975) p. 271.
\bibitem[KN98]{KN98} R. Kidambi and P. Newton, Motion of three point vortices on a sphere.
                     \emph{Physica D} {\bf 116} (1998), 143-175.
\bibitem[KN99]{KN99} R. Kidambi and P. Newton, Collapse of three vortices on a sphere.
                      \emph{Il Nuovo Cimento} {\bf 22} (1999), 779-791.
\bibitem[KN00]{KN00} R. Kidambi and P. Newton, Point vortex motion on a sphere with solid boundaries,
                      \emph{Phys. fluids} {\bf 12} (2000), 581-588.
\bibitem[KR89]{KR89} K. Klyatskin and G. Reznik, Point vortices on a rotating sphere,
                      \emph{Oceanology} {\bf 29} (1989), 12-16.
\bibitem[L65] {L65} H. Van Loon, A climatology study of the atmospheric circulation in the southern hemisphere
                      during the IGY, Part I: 1 July 1957-31 March 1958,
                      \emph{J. Appl. Meteo.} {\bf 4} (1965), 479-491.
\bibitem[L67] {L67} H. Van Loon, A climatology study of the atmospheric circulation in the southern hemisphere
                      during the IGY, Part II,
                      \emph{J. Appl. Meteo.} {\bf 6} (1967), 803-815.
\bibitem[LMR01]{LMR00} C. Lim, J. Montaldi, M. Roberts, Relative equilibria of point vortices on the sphere.
                        \emph{Physica D} {\bf 148} (2001), 97-135.
\bibitem[LMR]{LMR} F. Laurent-Polz, J. Montaldi, M. Roberts, Stability of point vortices on the sphere.
                       \emph{In preparation}.
\bibitem[LP02]{LP02} F. Laurent-Polz, Point vortices on the sphere: a case with opposite vorticities.
                       \emph{Nonlinearity} {\bf 15} (2002), 143-171.
\bibitem[LR96]{LR96} D. Lewis  and T. Ratiu, Rotating $n$-gon/$kn$-gon vortex configurations.
                        \emph{J. Nonlinear Sci.} {\bf 6} (1996), 385-414.
\bibitem[MPS99]{MPS99} J. Marsden, S. Pekarsky and S. Shkoller, Stability of relative equilibria of point vortices on a sphere and symplectic integrators.
                        \emph{Nuovo Cimento} {\bf 22} (1999).
\bibitem[M78]{M78} G. Mertz, Stability of body-centered polygonal configurations of ideal vortices.
                         \emph{Phys. Fluids} {\bf 21} (1978), 1092-1095.
\bibitem[Mi71]{Mi71} L. Michel, Points critiques des fonctions invariantes sur une $G$-vari\'{e}t\'{e}.
                         \emph{CR Acad. Sci. Paris} {\bf 272} (1971), 433-436.
\bibitem[MFG91] {MF91} C. Mechoso, J. Farrara, and M. Ghil. Intraseasonal variability of the
                       winter circulation in the southern hemisphere atmosphere.
                       \emph{J. Atmos. Sci.} {\bf 48} (1991), 1387-1404.
\bibitem[MR94] {MR94} J. Marsden and T. Ratiu, Introduction to Mechanics and Symmetry.
                       TAM {\bf 17} Springer-Verlag (1994).
\bibitem[MS71]{MS71} G. Morikawa and E. Swenson, Interacting motion of rectilinear geostrophic vortices.
                        \emph{Phys. Fluids} {\bf 14} (1971), 1058-1073.
\bibitem[N00] {N00} P. Newton, The $N$-vortex problem: analytical techniques.
                       {\bf 145} Springer-Verlag (2000).
\bibitem[Or98] {Or98} J-P. Ortega, Symmetry, Reduction, and Stability in Hamiltonian Systems.
                        Ph.D. Thesis. University of California, Santa Cruz (1998).
\bibitem[OR99]{OR99} J-P. Ortega, T. Ratiu, Stability of Hamiltonian relative equilibria.
                         \emph{Nonlinearity} {\bf 12} (1999), 693-720.
\bibitem[P79]{P79} R. Palais, Principle of symmetric criticality.
                         \emph{Comm. Math. Phys.} {\bf 69} (1979), 19-30.
\bibitem[Pa92] {Pa92} G. Patrick, Relative equilibria in Hamiltonian systems: the dynamic interpretation of nonlinear stability on a reduced phase space.
                       \emph{J. Geom. Phys.} {\bf 9} (1992), 111-119.
\bibitem[PM98]{PM98} S. Pekarsky and J. Marsden, Point vortices on a sphere: Stability of relative equilibria.
                         \emph{J. Math. Phys.} {\bf 39} (1998), 5894-5907.
\bibitem[PD93]{PD93} L. Polvani and D. Dritshel, Wave and vortex dynamics on the surface of a sphere.
                         \textit{J. Fluid Mech.} {\bf 255} (1993), 35-64
\bibitem[PRW02] {PRW01} G. Patrick, M. Roberts, and C. Wulff, Stability of Poisson equilibria and Hamiltonian relative equilibria by energy methods.
                       \emph{Preprint arXiv: math.DS/0201239} (2002).
\bibitem[RL53] {RL53} J. Rubin and H. Van Loon, Aspects of the circulation of the southern hemisphere.
                       \emph{J. Meteo.} {\bf 11} (1953), 68-76.
\bibitem[S43] {S43} H. Stewart, Periodic properties of the semi-permanent atmospheric pressure systems.
                       \emph{Quart. Appl. Math.} {\bf 1} (1943), 262-267.
\bibitem[S68] {S68} N. Streten, Some aspects of high latitude southern hemisphere summer circulation as
                       viewed by ESSA 3. \emph{J. Appl. Meteo.} {\bf 7} (1968), 324-332.
\end{thebibliography}
\end{document}